\documentclass{amsart}

\usepackage{mathrsfs}
\usepackage{amsmath}
\usepackage[all]{xy}
 
\usepackage{mathtools}
\usepackage{mathbbol}
\usepackage{wasysym}
\usepackage{amssymb}

\usepackage[textsize=small]{todonotes}

\usepackage{MnSymbol}
\usepackage{comment}
\usepackage{enumerate}
\usepackage{xcolor}

\usepackage{ushort}
 
\usepackage{changepage}
\usepackage{algorithm}
\makeatletter
\def\BState{\State\hskip-\ALG@thistlm}
\makeatother

\usepackage{pdflscape}
\usepackage[utf8]{inputenc}
\usepackage[textsize=small]{todonotes}

\usepackage{amsmath,latexsym,amssymb,amsthm,graphicx}
\usepackage{mathtools}
\usepackage{xfrac}
\usepackage{fix-cm}
\usepackage{color,transparent}
\usepackage{verbatim}
\usepackage[a-3u]{pdfx}

\usepackage{epic,bez123}
\usepackage{wrapfig}

\usepackage{soul}
\usepackage{url}
\usepackage{tikz}

\usepackage[margin=1.1in]{geometry}

 \author[Ding]{Lihuang Ding}	
 \email{shanquan2@stu.pku.edu.cn}
\address{Beijing International Center for Mathematical Research\\Peking University\\ Beijing 100871, China P.R.}

\author[Martínez-Granado]{Dídac Martínez-Granado}
\email{didac.martinezgranado@uni.lu}
\address{Department of Mathematics\\University of Luxembourg\\Av. de la Fonte 6, Esch-sur-Alzette, L-4364, Luxembourg}

\author[Zalloum]{Abdul Zalloum}
\email{abdul.zalloum@utoronto.ca}
\address{Department of Mathematics\\Queens University\\ 48 University Ave.\\Jeffery Hall\\
Kingston, ON Canada K7L 3N6}

\newtheorem{thm}{Theorem}

\newtheorem{mainthm}{Theorem}

\usepackage{hyperref} 
\usepackage{refcount}
\newenvironment{mainthmbis}[1]
  {
   \renewcommand{\themainthm}{\getrefnumber{#1}$'$}%
   \refstepcounter{mainthm}%
   \addtocounter{mainthm}{-1}%
   \phantomsection%
   \begin{mainthm}
  }
  {\end{mainthm}}

  \numberwithin{thm}{section}
\numberwithin{figure}{section}

\newtheorem{theorem}[thm]{Theorem}

\newenvironment{restatethm}[1]
  {
   \renewcommand{\themainthm}{\getrefnumber{#1}}%
   \refstepcounter{mainthm}%
   \addtocounter{mainthm}{-1}%
   \phantomsection%
   \begin{mainthm}
  }
  {\end{mainthm}}

\newtheorem{lem}[thm]{Lemma}
\newtheorem{claim}[thm]{Claim}
\newtheorem*{claim*}{Claim}
\newtheorem{cor}[thm]{Corollary}
\newtheorem{prop}[thm]{Proposition}

\newtheorem*{ex*}{Exercise}

\newtheorem*{fact*}{Fact}

\theoremstyle{definition}
\newtheorem{remark}[thm]{Remark}

\newtheorem{definition}[thm]{Definition}

\newtheorem{conjecture}[thm]{Conjecture}

\newtheorem{construction}[thm]{Construction}

\theoremstyle{definition}
\newtheorem{defn}[thm]{Definition}

\theoremstyle{remark}
\newtheorem{quest}[thm]{Question}

\newcommand{\CAT}{\mathrm{CAT}}

\newcommand{\isom}{\mathrm{Isom}}

\newcommand{\U}{X}

\newcommand{\hull}{\mathrm{hull}}

\newcommand{\f}{\mathscr F}

\newcommand{\ax}{\mathrm{Ax}}
\newcommand{\diam }[1]{{\textbf{diam}\big(#1\big)}}

\makeatletter
\let\c@equation\c@thm
\makeatother
\numberwithin{equation}{section}

\DeclareMathOperator{\Cay}{\mathit{Cay}}
\DeclareMathOperator{\MCG}{\mathit{MCG}}

\DeclareMathOperator{\DDG}{\mathcal{D}(G)}

\DeclareMathOperator{\DDGS}{\mathcal{D}_{s}(G)}

\DeclareMathOperator{\Isom}{Isom}

\usepackage{mathtools}
\usepackage{xfrac}
\usepackage{fix-cm}
\usepackage{color,transparent}
\usepackage{verbatim}
\usepackage[a-3u]{pdfx}

\newcommand{\DMGtodo}[2][]{\todo[color=orange!80,#1]{#2}}

\theoremstyle{plain}

\usetikzlibrary{arrows}
\tikzset{
    labl/.style={anchor=south, rotate=90, inner sep=.5mm}
}
\tikzstyle{every picture}=[> = to]
\tikzset{cdlabel/.style={execute at begin node=$\scriptstyle,execute at end node=$}}
\tikzset{implication/.style={double equal sign distance, -implies}}
\tikzset{biimplication/.style={double equal sign distance, implies-implies}}

\begin{document}

\title{Growth tightness and genericity for word metrics from injective spaces}


\begin{abstract} Mapping class groups are known to admit geometric (proper, cobounded) actions on injective spaces. Starting with such an action, and relying only on geometric arguments, we show that all finite generating sets resulting from taking large enough balls in the respective injective space yield word metrics where pseudo-Anosov maps are exponentially generic. We also show that growth tightness holds true for the Cayley graphs corresponding to these finite generating sets, providing a positive answer to a question by Arzhantseva, Cashen and Tao.
\end{abstract}
\maketitle



\section{Introduction}

\subsection{Motivating questions:} 
Given a proper action of a group $G$ on a metric space $(X,d)$ and a subset $A \subset G$, the \emph{critical exponent} of the $A$-action on $(X,d)$ is given by 

\begin{equation} \omega_A(d) \coloneqq \limsup_{k \rightarrow \infty}\frac{\log \#\{g \in A \;\big|\;  d(x_0,gx_0) \leq k\}}{k},
\label{eq:criticalexp}
\end{equation}
where $x_0 \in X$. We are interested in questions related to comparing $\omega_{A}(d)$ when $A$ varies in $G$:

\begin{enumerate}
\item \textbf{Growth gap:} Given a subset $A \subset G$ of distinguished elements with no ``hyperbolic-like" properties, do we have $\omega_{G}(d) > \omega_A(d)$?
\vspace{1mm}
\item \textbf{Growth tightness:} Given any $H \subset G$ infinite normal subgroup of $G$, there is an induced pseudometric $\bar{d}$ on $G/H$. Do we have $\omega_{G}(d) > \omega_{G/H}(\bar{d})$?
\vspace{1mm}

\item \textbf{Lower cogrowth:} Given any $H \subset G$ infinite normal subgroup of $G$, do we have $\omega_{H}(d) > \omega_{G}(d)/2$?
\end{enumerate}

Originally, these questions were posed in the setting where $(X, d)$ is a Cayley graph \cite{growth_tight}. In this context, many of these questions remained open for a long time, particularly when $G$ lies outside the realm of hyperbolic and relatively hyperbolic groups. For example, it was only in recent work by Choi~\cite{Choi2021PseudoAnosovsAE} that examples of Cayley graphs of mapping class groups were constructed that demonstrated certain desirable growth-gap properties.

On the other hand, recent work by multiple authors (e.g., \cite{Yang2018}, \cite{Yang17:ConvexCocompact}, \cite{Yan14:GrowthTight}, \cite{COULON2024}, \cite{CGZ24:Confined}, and \cite{Arzhantseva2015}) has shown that when a group $G$ admits a well-behaved action on a space $X$ that contains a strongly contracting direction, all of the questions above have affirmative answers when counting orbit points in a metric ball of $X$. For instance, in the case of the mapping class group's action on Teichmüller space.

One might hope that for every Cayley graph on the mapping class group, all pseudo-Anosovs are strongly contracting elements. However, this is not the case: in \cite{Rafi2021}, the authors provide an example of a finite generating set for the mapping class group of the five-times punctured sphere, where some pseudo-Anosov elements lack the strongly contracting property. 

Let $G$ be a group that acts properly and coboundedly on a space $X$ with a strongly contracting element. The primary goal of this article is to develop an approach that bridges the gap between our understanding of the growth of the $G$-orbit in $X$ and the growth of word metrics on $G$. This approach is outlined in Subsection \ref{subsection:approach}.

As a consequence, we construct the first known examples of Cayley graphs for the mapping class group that satisfy both exponential genericity and growth tightness properties. Moreover, our results apply to a much broader class of groups, encompassing mapping class groups.

\subsection{The space of all geometric actions of a given group}\label{subsection:approach}

Starting with a hyperbolic group $G$, Furman \cite{Fur02:Coarse} introduced the space $\DDG$, which consists of all left-invariant pseudo-metrics on $G$ that are quasi-isometric to a given word metric, considered up to rough similarity. Equivalently, $\DDG$ can be viewed as the space of geometric actions of $G$ on roughly geodesic metric spaces. This space was later studied by Reyes~\cite{OR22:SpaceMetric} and Cantrell–Reyes~\cite{COR22:Manhattan}, who initiated its study as a metric space. 

In the present article, we move beyond the setting of hyperbolic groups and consider the space of all proper, cobounded actions of an arbitrary finitely generated group $G$ on roughly geodesic metric spaces. That is, for a group $G,$ define
\[\DDG:=\{(X,d): (X,d) \text{ is roughly geodesic and }G\stackrel{Geom}{\curvearrowright}(X,d)
\}/ \sim,\]

where $(X,d) \sim (Y,d')$ if given $x_0 \in X$, $y_0 \in Y$ there exist positive constants $A,B$ such that $|A.d(x_0,gx_0)-d'(y_0,gy_0)|<B$ for all $g \in G$. In words, two metrics are equivalent if, after possibly scaling one of them, they differ by (at most) a uniform additive amount. Our primary goal is to study $\DDG$ when  $G$ is the mapping class group of a finite type surface, where $\DDG$ is known to contain especially interesting points, namely, various geometric actions on \emph{injective metric spaces} as shown by \cite{haettelhodapetyt:coarse} and \cite{PetytZalloum24} (see also Theorem \ref{thm:infinitely_many_injective}). Our interest in these points comes mainly from work of Sisto and Zalloum \cite{SistoZalloum2024} who showed that pseudo-Anosov mapping classes act as strongly contracting elements on these spaces (in fact, they are characterized by this property). In particular, while very little is known about counting with respect to word metrics for mapping class groups, as we mentioned above, a great deal is known in the presence of strongly contracting directions. Our approach is thus to use our knowledge of the growth of a mapping class group's action on an injective space $(X,d)$ and export it to eventual terms of a sequence $(X_i,d_i) \in \DDG$ that converges to $(X,d)$.

\begin{mainthm}\label{thm:Most_main_theorem} Let $G$ be a finitely generated group, $(X,d) \in \DDG$, and let $A_1,A_2 \subset G$ be so that $\omega_{A_1}(d)>c.\omega_{A_2}(d)$ for some $c>0$. There exists a sequence of finite generating sets $S_n$ whose corresponding word metrics $d_{S_n}$ satisfy $d_{S_n} \rightarrow d$, and for all $n$ large enough, we have
$\omega_{A_1}(d_{S_n})>c.\omega_{A_2}(d_{S_n}).$
\label{thm:main}
\end{mainthm}
 The space $\DDG$ will be equipped with a certain topology induced by  the stable length spectra, see Subsection \ref{subsec:metric} for more details. Below, we describe how to obtain the sequence $d_{S_n}$ in Theorem \ref{thm:Most_main_theorem} (and in most of the theorems in the present paper).

\begin{construction} Let $G$ be a finitely generated group and $(X,d) \in \DDG$ be a special point, for instance, an injective metric space, or more generally, a metric space with strongly contracting directions. One can always construct a sequence $d_{S_n}$ with $d_{S_n} \rightarrow d$ by simply taking larger and larger balls $B_X(x_0,n) \subset X$, letting $S_n=\{g \in G: g \cdot x_0 \in B_X(x_0,n)\}$ and applying Švarc-Milnor to obtain $\Cay(G,S_n)$ with the corresponding word metric $d_{S_n}$ (see \cite[Lemma~4.2.12]{Reyes23:Thesis}).
\label{cons:generatingsets}
\end{construction}

If one replaces the $\limsup$ by $\lim$ in Equation~\ref{eq:criticalexp}, the limit need not generally exist. However, there are natural situations where it does, for instance, when $(X,d) \in \DDG$ is a Cayley graph of a finite generating set, or when $(X,d) \in \DDG$ admits a strongly contracting element~\cite{Yang17:ConvexCocompact}. In such a case, we say that $\omega_A(d)$ is a \emph{true limit}. 

\begin{prop}\label{prop:rate_to_density} Let $G$ act properly on a geodesic metric space $(X,d)$ such that $\omega_G(d)$ is a true limit. For any $A \subset G$, we have $\omega_A(d)<\omega_G(d)$ if and only if the following converges exponentially quickly

\[ \frac{\#\{g \in G: g \notin A \text{ and }d(x_0,gx_0) \leq k\}}{\# \{g \in G: d(x_0,gx_0) \leq k\}} \rightarrow 1\] 
\label{prop:gapimpliesexpgenericity}
\end{prop}

We refer to Section~\ref{sec:proofs} for the proof of Proposition~\ref{prop:gapimpliesexpgenericity}.

\subsection{Exponential genericity of pseudo-Anosovs} When William Thurston first announced his construction of pseudo-Anosov mapping classes, he encountered significant resistance. In response, he asserted that not only do they exist, but they must also be generic. 

There are multiple ways to interpret what it means for pseudo-Anosov mapping classes to be generic. However, the most persistent interpretation—one that has remained open since Thurston's development of the theory of surface automorphisms—states the following:

\begin{conjecture}[{\cite[Conjecture~3.15]{FarbProblems}}]\label{conj:Farb} Let $\text{PA} \subset \MCG(\Sigma)$ denote the collection of pseudo-Anosov mapping classes, $S \subset G$ be a finite generating set and $e \in G$. As $k \rightarrow \infty$, we have 
\[ \frac{ |\text{PA} \cap B_S(e,k)|}{|B_S(e,k)|} \rightarrow 1.\]

\end{conjecture}

Now, if \( d_S \) denotes the word metric corresponding to a finite generating set \( S \) for \( G=\MCG(\Sigma) \), and if \( \text{Red} \subset G=\MCG(\Sigma) \) denotes the collection of reducible mapping classes (which are precisely the complement of pseudo-Anosovs in \( G \)), then the inequality  
$\omega_{\text{Red}}(d_S) < \omega_G(d_S)$
is strictly stronger than Conjecture \ref{conj:Farb} above. More precisely, if $\omega_{\text{Red}}(d_S) < \omega_G(d_S)$
, then the ratio of pseudo-Anosovs in a ball of radius \( n \) to the total size of that ball converges to 1 \emph{exponentially quickly} in $\Cay(G, S).$

Recall that mapping class groups (and more generally strong coarse median groups) admit proper cobounded actions on injective metric spaces (\cite{haettelhodapetyt:coarse}, \cite{PetytZalloum24}). In \cite{SistoZalloum2024}, and relying on work of Yang \cite{Yang2018}, Sisto and Zalloum proved that $\omega_{\text{Red}}({d_X}) < \omega_G(d_X)$ for any injective metric space $(X,d_X)$ on which $G=\MCG(\Sigma)$ acts geometrically. By exhibiting a sequence of word metrics for mapping class groups $d_{S_n} \in \DDG$ that converges to such an injective metric $(X,d_X)$, we obtain, from Theorem~\ref{thm:Most_main_theorem} and Proposition \ref{prop:rate_to_density}, the following result in the setting of $\MCG(\Sigma)$.

\begin{mainthm}\label{thm:p.A are dense} Let $X$ be an injective space on which $\MCG(\Sigma)$ acts properly coboundedly. There exists $n_0$ such that for each $n \geq n_0$, the finite set $S_n=\{g \in G: g \cdot x_0 \in B_X(x_0, n) \subset X\}$ generates $\MCG(\Sigma)$ and satisfies $\omega_{\text{Red}}(d_{S_n}) < \omega_G(d_{S_n}).$ In particular, Conjecture~\ref{conj:Farb} holds true in $\Cay(\MCG(\Sigma), S_n)$ with the convergence happening exponentially quickly.
\end{mainthm}

The first major breakthrough regarding Conjecture \ref{conj:Farb} comes from Choi's work in \cite{Choi2021PseudoAnosovsAE}. Starting with an arbitrary finite generating set \(S\), Choi enlarges \(S\) to a finite set \(S'\) by adding large powers of a fixed pair of pseudo-Anosov elements, allowing him to show that $\omega_{\text{Red}}({d_{S'}}) < \omega_G({d_{S'}})$. 

The finite generating sets \( S_n \) constructed in Theorem \ref{thm:p.A are dense} are significantly more natural than those in Choi's work \cite{Choi2021PseudoAnosovsAE}, where his generating sets consist primarily of large powers of two pseudo-Anosov elements. In contrast, the finite generating sets in Theorem \ref{thm:p.A are dense} arise simply by taking sufficiently large balls in \( X \) and applying the Švarc–Milnor lemma. 

Two months after we posted the first draft of this paper to the arXiv, Choi fully proved Conjecture \ref{conj:Farb} using very different methods \cite{Full}. Specifically, he showed that for \emph{every} finite generating set \( S \) of \( G=\MCG(\Sigma) \), the ratio of pseudo-Anosov elements in a ball of radius \( n \) in \( \Cay(G,S) \) to the total number of elements in that ball approaches 1 as \( n \to \infty \). However, his result does not establish that $\omega_{\text{Red}}({d_S}) < \omega_G({d_S})$, and in particular, it does not imply exponential convergence. 

Thus, if one adopts the strictly stronger version of Conjecture \ref{conj:Farb}—namely, that for every finite generating set \( S \) of \( G=\MCG(\Sigma) \), one must have \( \omega_{\text{Red}}({d_S}) < \omega_G({d_S}) \)—then the only known results in the literature remain Theorem \ref{thm:p.A are dense} and Choi's work in \cite{Choi2021PseudoAnosovsAE}. As discussed earlier, Theorem \ref{thm:p.A are dense} has the advantage that its finite generating sets are significantly more natural than those provided in \cite{Choi2021PseudoAnosovsAE}.

\begin{quest} 
Do we have \( \omega_{\text{Red}}({d_S}) < \omega_G({d_S}) \) for \emph{every} finite generating set $S$ of $G=\MCG(\Sigma)$?
\end{quest}

It is worth noting that while pseudo-Anosov mapping classes act on Teichmüller space as strongly contracting elements (via work of Minsky \cite{Minsky1996} for the Teichmüller metric, Lenzhen, Rafi and Tao \cite{Lenzhen2011} in the Thurston metric), the action is not cobounded. Therefore, we can't apply our strategy using Teichmüller space instead of the injective space $X$ above.

\begin{cor} For any finite $A \subset \MCG(\Sigma)$, there exists a finite set $S$ with $A \subset S \subset \MCG(\Sigma)$ so that pseudo-Anosov elements are exponentially generic in $\Cay(\MCG(\Sigma), S)$. The generating set $S$ is constructed by picking an arbitrary point $x_0$ in an arbitrary injective dual of $\MCG(\Sigma)$, taking large enough $n >0$ so that $A. x_0 \subset B(x_0,n)$ and $S \coloneqq \{ g \in G | g . x_0 \in B(x_0,n) \}$ generates $G$.
\end{cor}

It is also worth noting that pseudo-Anosov mapping classes were proven by Maher to be generic with respect to random walks, see \cite{maher:random} for more details.

\subsection{Injective metric spaces}\label{subsec:construction_injective}

In this section, we provide some preliminaries on injective metric spaces and discuss the specific injective spaces we will be working with.

\begin{definition} A metric space $X$ is \emph{injective} if it is the injective object in the category of metric spaces with respect to 1-Lipshitz maps. That is to say, $X$ is injective if for any 1-Lipshitz map $f:Y \rightarrow X$ and any isometric embedding $i:Y \hookrightarrow Z$; there is a $1$-Lipshitz map $\tilde{f}: Z \rightarrow X$ such that $\tilde{f}\circ i=f.$
\end{definition}

Observe that injective metric spaces are geodesic. Namely, if $X$ is injective and $x,y \in X$, one can consider the two-point metric space $Y=\{x,y\}$ and the interval $Z=[0, d(x,y)]$. We then have the isometric embeddings $f:Y \rightarrow X$ and $i: Y \rightarrow Z$ and using the definition of injective spaces, we get a map $\tilde{f} \colon [0, d(x,y)] \rightarrow X$ which is 1-Lipschitz and is an isometry on the end points with $\tilde{f}(0)=x,\, \, \tilde{f}(d(x,y))=y$. The triangle inequality then implies $\tilde{f}$ is an isometric embedding (cf.\cite[Page~6]{Lang13:Injective}).
There is another equivalent---and perhaps more hands-on---definition of injective spaces that we will give in Appendix~\ref{subsec:coarsely_injective}, and it is the one used in the following result. The two definitions are equivalent by \cite[Proposition~2.3]{Lang13:Injective}.

The main theorem we will use regarding injective spaces is the following:

\begin{theorem}[Lemma 3.10 in \cite{haettelhodapetyt:coarse}] Every hierarchically hyperbolic group admits a proper cobounded action on an injective metric space.
\end{theorem}

To understand the types of mapping classes that appear in a ball of radius $n$ in Theorem \ref{thm:p.A are dense}, leading to exponential density, it is crucial to examine the origin of the injective metric space $X$. It is worth noting, however, that there exist various distinct injective metric spaces on which $\MCG(\Sigma)$ acts properly coboundedly (see Theorem \ref{thm:infinitely_many_injective}). Nevertheless, there is a particularly canonical one that can be easily described. See Appendix~\ref{subsec:proof_infinitely} as well as~\cite{PetytZalloum24} for more details.

For $x, y \in \MCG(\Sigma)$, the \emph{distance formula} (\cite{masurminsky:geometry:1}, \cite{masurminsky:geometry:2}) states that $d(x,y)$ is coarsely the sum of their large projections onto certain hyperbolic spaces, specifically the curve graphs of subsurfaces. There are finitely many such subsurfaces, denoted $Y_1, \dots, Y_n$. If $\{Y_{n_1}, \dots, Y_{n_k}\} \subset \{Y_1, \dots, Y_n\}$ is a maximal collection of curve graphs associated with pairwise disjoint subsurfaces, we modify the sum of distances by replacing $\sum_{1 \leq i \leq k} d_{Y_{n_i}}(x,y)$ with $\max_{1 \leq i \leq k} d_{Y_{n_i}}(x,y)$. The new total sum is strictly smaller than the original, but since there is a uniform bound on the number of pairwise disjoint subsurfaces, the original sum remains quasi-isometric to the modified one. 

The primary \emph{coarsely injective} metric space on which the mapping class group acts properly coboundedly is then obtained (roughly speaking) by defining $d(x,y)$ to be the aforementioned modified sum. While this only provides a coarsely injective space, such spaces are sufficient to establish the full conclusion of Theorem \ref{thm:p.A are dense}. Indeed, the only property of injective spaces used in the proof of Theorem \ref{thm:p.A are dense} is the existence of strongly contracting directions, a property that remains valid in the coarsely injective space we just described. 

Furthermore, every coarsely injective space admits an isometric embedding into a certain injective space—its \emph{injective hull}—with a coarsely surjective image. Consequently, the notions of injective and coarsely injective spaces are largely interchangeable. See Appendix~\ref{subsec:proof_infinitely} for more details.

\subsection{Growth Tightness} 

Growth tightness is a condition comparing the critical exponent of a given group with the critical exponent of its quotients. It was originally introduced by de la Harpe and Grigorchuk in the context of word metrics as a strategy to show that a given group is Hopfian~\cite{growth_tight}. 

The first general result appeared in Arzhantseva and Lysenok~\cite{AL02:GrowthTight}, who established that all word metrics in non-elementary hyperbolic groups are growth tight. Sambusetti~\cite{Sam02:GrowthTight_Surface} proposed extending growth tightness beyond word metrics to more general left-invariant pseudo-metrics, an idea later echoed by Yang~\cite{Yan14:GrowthTight}; Arzhantseva, Cashen, and Tao~\cite{Arzhantseva2015}, following~\cite{DPPS11:GrowthTight}. 

Let $G$ be a finitely generated group, and let $(X,d)$ be a space on which $G$ acts properly and coboundedly. The metric $d$ induces a pseudometric $\bar{d}$ on any quotient $G/H \cdot x_0$ of $G$, defined by  
\[
\bar{d}(g H \cdot x_0, g' H \cdot x_0) \coloneqq \inf_{h, h' \in H} d(g h\cdot x_0, g' h' \cdot x_0).
\]
Similarly, let $\omega_{G/H}(\bar{d})$ denote the growth exponent of $G/H \cdot x_0$. We say that $G$ is \emph{growth tight} with respect to $d$ if 
\[
\omega_{G/H}(\bar{d}) < \omega_{G}(d) 
\] 
for every infinite normal subgroup $H \triangleleft G$. 

Particularly relevant are metrics with negatively curved directions, where contracting or strongly contracting elements play a key role. Growth tightness for co-compact Kleinian groups was established by Sambusetti~\cite{Sam08:GrowthTight_Cocompact} and later extended to geometrically finite Kleinian groups with the \emph{parabolic gap property} by Dal'Bo, Peigné, Picaud, and Sambusetti~\cite{DPPS11:GrowthTight}. This property was further generalized by Yang~\cite{Yan14:GrowthTight} for cocompact actions with contracting elements, and Arzhantseva, Cashen, and Tao~\cite{Arzhantseva2015} for (posibly non-cocompact) actions with contracting elements and complementary growth gap, leading to the concept of statistically convex-cocompact actions formalized in~\cite{Yang17:ConvexCocompact} and further explored in~\cite{DY24:GrowthTight}. In particular, Arzhantseva, Cashen, and Tao~\cite{Arzhantseva2015} proved that the action of $\MCG$ on Teichmüller space is growth tight.

Growth tightness is a subtle and delicate property. For instance, Dal'Bo, Otal, and Peigné~\cite{DOP00:GrowthTight_Counter} constructed examples of groups acting on two different equivariantly quasi-isometric spaces, where one action is growth tight and the other is not. Similarly, in the setting of word metrics, Cashen and Tao~\cite{CT16:GrowthTight_Word} constructed examples of groups that exhibit growth tightness with respect to some generating sets but not others.

In~\cite[Question~3, Page 33]{Arzhantseva2015}, the authors posed the question of whether $\MCG(\Sigma)$ admits any finite generating set $S$ such that $\MCG(\Sigma)$ is growth tight with respect to $d_S$. 

Despite significant progress in understanding growth tightness in recent years, this question has remained open until now.

In this paper, we provide an affirmative answer to the question of Arzhantseva, Cashen, and Tao by proving that if a word metric is sufficiently close to an action with strongly contracting elements, then such a word metric must be growth tight.

In fact, we construct infinitely many generating sets for which growth tightness holds for a broader class of subgroups of $G$, known as \emph{confined subgroups}. 

Confined subgroups were originally introduced by Hartley and Zalesski\u{i}~\cite{HZ97:Confined} and have been extensively studied in recent works, including~\cite{LM18:Confined2, LM22:Confined3, FG23:ConfinedInfiniteVolume, GL24:ConfinedStationary, CGZ24:Confined}.

\begin{defn}
	Let $G$ be a locally compact, metrizable topological group, and $\Gamma \leq G$ a subgroup. A subgroup $H \leq G$ is called \textit{confined} by $\Gamma$ in $G$, with a finite \textit{confining subset} $P \subset G$ if, for every element $g \in \Gamma$,  
	\[
	g^{-1} H g \cap P \setminus \{1\} \neq \emptyset.
	\]
\end{defn}

In our setting, we will be working with $G$ discrete, and $\Gamma=G$, in which case we just say $H$ is a \emph{confined subgroup} of $G$, with confining subset $P$.
Note that the direct product of a finite normal subgroup with any subgroup is a confined subgroup. Let $E(G)$ denote the maximal finite normal subgroup of $G$. To eliminate such pathological examples, we focus on confined subgroups with \emph{non-degenerate} confining subsets $P \subset G$, where $P$ is disjoint from $E(G)$. Examples of groups with non-degenerate confining subsets include normal subgroups of infinite index.

\begin{mainthm}\label{GrowthTight4SomeS}
Let $X$ be an injective dual of a mapping class group $G$ of a finite-type surface. Then there exists $n_0$ such that for each $n \geq n_0$, the finite generating set $S_n=\{g \in G: g \cdot x_0 \in B_X(x_0,n)\}$ generates $G$ and the word metric $d_{S_n}$ satisfies
\[
\omega_{G/H}(d_{\bar{S}}) < \omega_G(d_S)
\] for any infinite confined subgroup $H\subset G$ with non-degenerate confining subset. 
\end{mainthm}

\subsection{Lower cogrowth bound}

Let $H$ be a confined subgroup of $G$ with a non-degenerate confining subset. The ratio $\omega_{H}(d_S)/\omega_{G}(d_S)$ is called the \emph{cogrowth} of $Q \coloneqq G / H$.

When $H$ is an infinite-index normal subgroup and $X = G = F_n$, the strict lower bound on cogrowth was established by Grigorchuk~\cite{Gri77:Symmetric}. Grigorchuk and de la Harpe~\cite[Page~69]{growth_tight} posed the question of whether this strict lower bound on cogrowth also holds when $F_n$ is replaced by a non-elementary Gromov hyperbolic group and $X$ is one of its Cayley graphs. This long-standing problem was recently resolved affirmatively by Jaerisch, Matsuzaki, and Yabuki~\cite{JMY19:Normalizer}.

In Theorem~\ref{thm:cogrowth_mcg} below, relying on recent work in~\cite{CGZ24:Confined}, we establish the existence of many Cayley graphs of mapping class groups that satisfy certain lower cogrowth inequalities. We state the result here for mapping class groups, with the more general version appearing in Subsection~\ref{subsec:generalities}.

\begin{mainthm}\label{thm:cogrowth_mcg}
Let $G=\MCG(\Sigma)$, for $\Sigma$ a finite type hyperbolic surface. If $H$ is a confined subgroup with non-degenerate confining subset,
then there exist generating sets $S_n$ of $G$ and $n_0>0$, so that for all $n \geq n_0$, 
then $\omega_{H}(d_{S_n}) > \omega_{G}(d_{S_n})/2.$
\end{mainthm}

Note that this result does not settle the cogrowth gap for these generating sets, since the inequalities are not uniform in $H$, but we expect techniques similar to those in the proof of Theorem~\ref{GrowthTight4SomeS} to work in this setting. This will be explored in a forthcoming paper.

\subsection{The space of geometric actions with strongly contracting elements}\label{subsec:metric} Let $G$ be a finitely generated group. Recall that $\DDG$ is defined as the space of all proper cobounded actions on roughly geodesic metric spaces. We will now introduce a topology on $\DDG$ induced by a (pseudo-)distance. This topology will be mainly used to motivate Question~\ref{que:densityinjective} below, while the distance on $\DDG$ will play a role in Theorem~\ref{thm:infinitely_many_injective}. We define an equivalence relation $\sim$ on $\DDG$ by setting $(X_1, d_1) \sim (X_2, d_2)$ if there exists a constant $c$ such that 
\[
\ell_{d_1}(g) = c \cdot \ell_{d_2}(g) \quad \text{for all } g \in G,
\]
where $\ell_d(g) \coloneqq \lim_{n \to \infty} d(x, g^n x)/n$ denotes the \emph{stable length} of $g$. Intuitively, two metrics are identified if they have the same projectivized length spectrum.

For $d, d' \in \DDG / \sim$, motivated by~\cite{Fur02:Coarse,OR22:SpaceMetric}, we define a metric on $\DDG / \sim$ as follows:

\[
D(d, d') := \log \sup_{g \in G_0} \frac{\ell_d(g)}{\ell_{d'}(g)} + \log \sup_{g \in G_0} \frac{\ell_{d'}(g)}{\ell_d(g)},
\]

where $G_0$ denotes the collection of all elements in $G$ with positive translation length in both (or equivalently, any) of the metrics $d$ and $d'$. It is straightforward to verify that $D$ defines a metric.

We topologize $\DDG$ with the coarsest topology such that the quotient map $q: \DDG \to \DDG / \sim$ is continuous. Let $\DDGS$ denote the subset of points $(X, d) \in \DDG$ such that $X$ contains an infinite strongly contracting rough geodesic. In light of length spectrum rigidity for actions with strongly contracting elements~\cite{W24:Marked}, the restriction of $q$ to $\DDGS$ is a bijection onto its image. Consequently, the subspace $\DDGS \subset \DDG$ inherits the same metric structure as defined above. However, since length spectrum rigidity is not known for all points in $\DDG$, the metric $D$ does not (a priori) extend to a metric on all of $\DDG$.

One can also use properties of a converging sequence to deduce information about its limit point. For example, in the setting of Gromov hyperbolic groups, in forthcoming work~\cite{CMGR:GreenDensity}, the second author and collaborators show that the space of all Green metrics is dense in $\DDG$. They use this result to establish certain dynamical properties for all points in $\DDG$ by leveraging the strong dynamical properties of the class of Green metrics.

In summary, the density of a certain subclass of metrics $Y \subset \DDG$ is a powerful tool both for transferring information from a limit point to its converging sequence and for extending properties of $Y$ to the entire space of metrics $\DDG$. This motivates the following question:

\begin{quest}
For the mapping class group $G = \MCG(\Sigma)$, is the subspace $\DDGS$ dense in the space of all metrics $\DDG$?
\label{que:densityinjective}
\end{quest}

A positive answer to Question~\ref{que:densityinjective} would imply the density of the approximating word metrics $S_n$ in Theorem~\ref{thm:genericity}, thereby yielding a dense subset of word metrics exhibiting both exponential genericity and growth tightness. Indeed, one of the crucial observations of our argument will be that for every $d_c \in \DDGS$, there exist a sequence of word-metrics $d_{S_n}$, with $d_{S_n} \to d_c$ in $\DDG$ (see Lemma~\ref{lem:density}) and, for $n$ large enough, all $d_{S_n}$ inherit the same desirable growth properties as $d_c$. Then, by density of $\DDGS$, a double diagonal argument would allow us to approximate any $d \in \DDG$ by a sequence of such word-metrics.

The following theorem demonstrates that the space $(\DDGS, D)$ is indeed large, at least when $G=\MCG(\Sigma)$:

\begin{mainthm}\label{thm:infinitely_many_injective}
Let $G=\MCG(\Sigma)$.
There exists a sequence of
injective metric spaces $(X, \sigma_k) \in \DDGS$ on which $G$ acts properly and coboundedly such that 
\[
D(\sigma_1, \sigma_k) \to \infty \quad \text{as } k \to \infty.
\]
In particular, 
\[
\diam \DDGS = \infty.
\]
\end{mainthm}

As Mahan Mj likes to say, there are two standard ways to ``kill" a hyperbolic surface (i.e., to go to infinity in Teichmüller space): either by \emph{throttling} the surface (letting the hyperbolic/translation length of a curve go to zero) or by \emph{wringing its neck} (applying large powers of a Dehn twist). We prove Theorem~\ref{thm:infinitely_many_injective} by throttling the curve graph of a subsurface (along with its translates). To achieve this, the machinery of \emph{curtains} (see~\cite{PSZCAT}, \cite{zalloum:injectivity}, and \cite{PetytZalloum24}) proves to be quite useful.

The key idea is as follows: nearest-point projections to geodesics in a given curve graph provide natural tripartitions of $\MCG(\Sigma)$. Specifically, if $b \subset \mathcal{C}(Y)$ is a geodesic and $I_k \subset b$ is a subgeodesic of length $k$, we obtain a tripartition of $\mathcal{C}(Y)$ based on whether a point $x \in \mathcal{C}(Y)$ projects inside, to the right, or to the left of $I_k$. This induces a collection of tripartitions in $\MCG(\Sigma)$ via the pullback under $\pi_Y: \MCG(\Sigma) \to \mathcal{C}(Y)$. The set of points in $\MCG(\Sigma)$ mapping to $I_k$ is called a \emph{$k$-curtain}.

Now, fix a subsurface $Z_1 \subseteq \Sigma$ and fix $h \in \MCG(\Sigma)$ be a partial pseudo-Anosov on $Z_1,$ in particular, $h$ is loxodromic on $\mathcal{C}(Z_1).$ Let $\mathcal{Z} = \{\mathcal{C}(Z_1), \dots, \mathcal{C}(Z_m)\}$ be the set of curve graphs where $g$ is loxodromic (the set is finite as the subsurfaces $Z_i$ need to be pairwise disjoint in this case) and consider the set

\[
\MCG(\Sigma) \cdot \mathcal{Z} = \{g \cdot \mathcal{C}(Z_i) \mid g \in \MCG(\Sigma), \mathcal{C}(Z_i) \in \mathcal{Z} \}.
\]

As discussed above, curtains in $\MCG(\Sigma)$ are generally defined by considering geodesics in curve graphs $\mathcal{C}(Y)$. We will now define a certain collection of curtains in $\MCG(\Sigma)$ of varrying ``thickness" depending on whether their defining curve graph lies in $\MCG(\Sigma) \cdot \mathcal{Z}$. More precisely, for $\mathcal{C}(Y) \in \MCG(\Sigma) \cdot \mathcal{Z}$; we consider the collection of all curtain in $\MCG(\Sigma)$ defined by subgeodesic $I_k$ of length $k$ in $\mathcal{C}(Y).$ Lets denote this collection by $\mathcal{H}^k$. On the other hand, for $\mathcal{C}(Y) \notin \MCG(\Sigma) \cdot \mathcal{Z}$, we consider the collection of all curtains in $\MCG(\Sigma)$ defined by subgeodesic $I_1$ of length $1$ in $\mathcal{C}(Y)$ and we let $\mathcal{H}^1$ denote such a collection. Finally, let $\mathcal{W}^k:= \mathcal{H}^1 \cup \mathcal{H}^k$.

We define a metric $\sigma_k$ on $\MCG(\Sigma)$ by counting the maximal chain of curtains in $\mathcal{W}^k$ separating two points. By the work of~\cite{PetytZalloum24}, one can verify that the corresponding metric space, denoted $(\MCG(\Sigma), \sigma_k)$, is coarsely injective. The crucial observation is that, by definition of $\mathcal{Z}$, every curve graph where $h$ is loxodromic must belong to $\mathcal{Z}$. However, since the curtains arising from $\mathcal{Z}$ are at least $k$-thick, the translation length of $h$ in $(\MCG(\Sigma), \sigma_k)$ approaches zero as $k \to \infty$. See Appendix~\ref{subsec:proof_infinitely} for further details.

\begin{remark}
To carry out the argument above, it is not necessary for $Z_1$ to be a proper subsurface. One can simply take $Z_1 := \Sigma$ and choose $h$ to be a pseudo-Anosov homeomorphism on $\Sigma$. This choice simplifies the computation, as the set $\mathcal{Z}$ becomes a singleton, and so does $\MCG(\Sigma) \cdot \mathcal{Z}$. Accordingly, in the proof provided in Appendix~\ref{subsec:proof_infinitely}, we adopt this simplified setup. We chose to present the more general case above, however, as it may prove useful in other contexts.
\end{remark}

\subsection{Arbitrarily good hierarchy paths} 

A well-known open problem asks whether the mapping class group of a finite-type surface admits a finite generating set such that every pair of points is connected by a geodesic (or at least a rough geodesic) that is also a hierarchy path. 

On the other hand, in recent work by the third author with Petyt~\cite{PetytZalloum24}, a coarsely injective metric $\sigma$ on $\MCG(\Sigma)$ was constructed, where every pair of points $g, h \in \MCG(\Sigma)$ is connected by a hierarchy path that is also an $\alpha$-rough geodesic for some uniform constant $\alpha$. 

In the present article, by exhibiting a sequence of word metrics $d_{S_n} \to \sigma$, we establish the following result:

\begin{mainthm}\label{thm:good_hierarchies} 
For each $\epsilon > 0$, there exists a finite generating set $S_n$ for $\MCG(\Sigma)$ such that every pair of points $g,h \in \MCG(\Sigma)$ is connected by a path $\gamma_{gh}$ satisfying:

\begin{enumerate}
    \item $\gamma_{gh}$ is a hierarchy path, and
    
    \item $\gamma_{gh}$ is a $(1+\epsilon, 1)$-quasi-geodesic in $\Cay(\MCG(\Sigma), S_n)$.
\end{enumerate}

\end{mainthm}

\subsection{Generalities}
\label{subsec:generalities}
We first stated our results in the setting of $G = \MCG(\Sigma)$ for a finite-type surface $\Sigma$, as this setting motivates most of the original questions. However, our results hold in a much more general framework where $G$ acts properly coboundedly on a metric space in which all Morse geodesics are (quantitatively) strongly contracting. For properties of such spaces, see \cite{Yang2018}, \cite{Yang17:ConvexCocompact}, \cite{COULON2024}, and \cite{CGZ24:Confined}. Spaces satisfying this property include CAT(0) spaces~\cite{charneysultan:contracting}, injective spaces~\cite{SistoZalloum2024}, Garside groups~\cite{calvezwiest:morse}, certain small-cancellation groups~\cite{zbinden:small}, and weakly modular graphs with convex balls~\cite{soergelzalloum:morse}.

Our Theorem~\ref{thm:p.A are dense} follows from the following more general result (proven in Section~\ref{sec:proofs}) combined with~\cite{Yang2018}. 

\begin{mainthmbis}{thm:p.A are dense} \label{thm:genericity} 
Let $G$ act properly coboundedly on a geodesic space $(X,d)$ with a strongly contracting element, and let $x_0 \in X$. Define $A \subset G$ to be the set of elements that are \emph{not} strongly contracting on $X$. Then, there exists $n_0 > 0$ such that for all $n \geq n_0$, if $S_n$ consists of elements whose $x_0$-orbit lies in $B(x_0,n) \subset X$, then $S_n$ generates $G$ and satisfies 
\[
\omega_{A}(d_{S_n}) < \omega_G(d_{S_n}).
\]
Moreover, as $k \to \infty$, the following convergence occurs exponentially quickly:
\[
\frac{\#\{g \in G : g \text{ is strongly contracting and } d_{S_n}(e,g) \leq k\}}{\# \{g \in G : d_{S_n}(e,g) \leq k\}} \to 1.
\]
\end{mainthmbis}

Similarly, our growth tightness result holds in a much broader context.

\begin{mainthmbis}{GrowthTight4SomeS}\label{mainthm:cogrowth}
Let $G$ act properly coboundedly on a geodesic space $(X,d)$ with a strongly contracting element, and let $x_0 \in X$. Then, there exists $n_0 > 0$ such that for all $n \geq n_0$, if $S_n$ consists of elements whose $x_0$-orbit lies in $B(x_0,n) \subset X$, then $S_n$ generates $G$ and satisfies 
\[
\omega_{G/H}(d_{\bar{S}}) < \omega_G(d_S)
\]
for any infinite confined subgroup $H \subset G$ with a non-degenerate confining subset.
\end{mainthmbis}

The proof of Theorem~\ref{mainthm:cogrowth} is given in Section~\ref{sec:proof_GrowthTightness}.

\subsection*{Organization}
Section~\ref{sec:proofs} (pages 8-10) contains the proofs of exponential genericity (Theorems~\ref{thm:main},~\ref{thm:p.A are dense}, and Theorems~\ref{thm:genericity}), hierarchies (Theorem~\ref{thm:good_hierarchies}) and cogrowth (Theorem~\ref{thm:cogrowth_mcg}). Section~\ref{sec:proof_GrowthTightness} (pages 10-19) presents the argument for growth tightness, proving Theorems~\ref{GrowthTight4SomeS} ,~\ref{mainthm:cogrowth}.
Finally, Appendix~\ref{subsec:proof_infinitely} (pages 20–23) contains the proof of the infinite diameter result (Theorem~\ref{thm:infinitely_many_injective}).

\section*{Acknowledgements}

The authors are deeply grateful to Inhyeok Choi for his valuable comments and suggestions on earlier versions of the paper, as well as for identifying a gap in the proof of Theorem \ref{GrowthTight4SomeS} in the first draft.
DMG and AZ are thankful to Eduardo Reyes for many fruitful discussions, and especially for pointing out Lemma 4.2.12 of his PhD thesis where he keeps track of the constants in the Švarc-Milnor lemma. They are also quite grateful to Matthew Durham, Harry Petyt, Sam Taylor, Giulio Tiozzo and Wenyuan Yang for some productive feedback on the first draft of the article. They also wish to also thank Ilya Gekhtman, Christian Gorski and Kasra Rafi for interesting questions and discussions. LD is greatly thankful to Wenyuan Yang for fruitful discussions and feedback on the article. Finally, they would like to express their deep gratitude to the Fields Institute for their hospitality during the Thematic Program on Randomness and Geometry where this work was carried out. DMG was supported by the Luxembourg National Research Fund (FNR), grant reference AFR/Bilateral-ReSurface 22/17145118, and by the Long Term Visitor fund of the Thematic Program on Randomness and Geometry at the Fields Institute.
Finally, we thank the referee for their critical reading of the first version of the paper, and for also spotting the aforementioned gap in our previous proof of Theorem \ref{GrowthTight4SomeS} in the first draft.

\section{Proof of exponential genericity}
\label{sec:proofs}

In this section we give the proof of the exponential genericity result. We state the well-known Švarc-Milnor lemma (for a precise reference of the lemma with the constants below, see~\cite[Lemma~4.2.12]{Reyes23:Thesis}). The statement in ~\cite[Lemma~4.2.12]{Reyes23:Thesis} is given in the context of hyperbolic groups, but the proof doesn't use that assumption.

\begin{lem}\label{lem:density} Let $(X,d)$ be an $\alpha$-roughly geodesic metric space on which $G$ acts properly coboundedly and let $x_0 \in X$. If $m$ is large enough so that $S_m \coloneqq \{ g \in G | g. x_0 \in B(x_0,m) \}$ generates $G$, then, for any $ n \geq m$, we have $$(n-1-\alpha)d_{S_n}(g,h)-(n-1) \leq d(gx_0,hx_0) \leq n d_{S_n}(g,h).$$
\end{lem}

The following is an immediate consequence:

\begin{cor}\label{cor:convergence} Let $G$ act properly co-boundedly on some $\alpha$-roughly geodesic metric space $(X,d)$, the sequence of word metrics in the previous lemma $d_{S_n}$ converges to $d$ in $\DDG$. In fact, for any $\epsilon>0,$ there exists $n_0$ such that for all $n \geq n_0,$ we have the following (up to scaling the metrics) 
$$1-\epsilon \leq \frac{\ell_{d_{S_n}}(g)}{\ell_{d}(g)} \leq 1$$
for every $g \in G_0$.
    
\end{cor}

\begin{remark} One might also be interested in the critical exponent corresponding to counting conjugacy classes of a group $G$ that acts properly coboundedly on a certain metric space $(X,d)$. For a subset $A \subset G$, a point $x_0 \in X$ define

\[
\omega_A^c(d)\coloneqq 
\limsup_{k\to\infty}
\frac{\log \#\{\, [g]\in A \;\big|\; \ell_d^{x_0}([g]) \le k \,\}}{k}.
\]
where $\ell_d^{x_0}([g])= \text{min}\{d(x_0, g \cdot x_0): g \in [g]\}$. Work of \cite{Gekhtman2022} shows that if $G$ admits a proper cobounded action on a space $X$ with a strongly contracting element, and $A \subset G$ is the collection of elements that are not strongly contracting on $(X,d)$, then $\omega_A^c(d)< \omega_G^c(d)$. A straight forward computation using Corollary \ref{cor:convergence} shows that there are infinitely many finite generating sets $S_n$ such that $\omega_A^c(d_{S_n})< \omega_G^c(d_{S_n})$.
\end{remark}

The following is also immediate by Lemma \ref{lem:density}.

\begin{lem} Continuing with the notation of the Lemma~\ref{lem:density}, given $n \geq m$, we have 

$$(n-2-\alpha)d_{S_n}(g,h) \leq d(gx_0,hx_0) \leq n d_{S_n}(g,h),$$ for all $g,h$ with $d_{S_n}(g,h) \geq n-1. $
\label{lem:cutoff}
\end{lem}

\begin{proof}
    Since $g,h$ satisfy $d_{S_n}(g,h) \geq n-1$, the Lemma follows immediately by Lemma~\ref{lem:density} .
\end{proof}

\begin{lem} Continuing with the notation from Lemma~\ref{lem:density}, for a subset $A \subset G,$ let $\omega_A (d), \omega_A(d_{S_n})$ be the critical exponents of $A.x_0 \subset X, A.e \subset \Cay (G,S_n)$ (for $e \in G$) respectively. We have $$\frac{\omega_A(d_{S_n})}{n} \leq \omega_A (d) \leq \frac{\omega_A(d_{S_n})}{n-2-\alpha}$$
for every $n \in \mathbb{N}$. 

    \label{lem:criticalexponents}
\end{lem}

\begin{proof} By Lemma~\ref{lem:cutoff}, for $n \geq m$, we have $(n-2-\alpha)d_{S_n}(g,h) \leq d_X(gx_0,hx_0) \leq n d_{S_n}(g,h),$ for all $g,h$ with $d_{S_n}(g,h) \geq n-1$. Since the critical exponent is not affected by removing elements in a fixed ball, for each $n$, we have $\frac{\omega_A(d_{S_n})}{n} \leq \omega_A (d) \leq \frac{\omega_A(d_{S_n})}{n-2-\alpha}.$
\end{proof}

\subsection{Exporting growth rates to converging sequences} We now give a proof to Theorem \ref{thm:main}

\begin{proof} Let $A_1,A_2$ be the subsets of the statement. By Lemma~\ref{lem:criticalexponents}, for $i \in \{1,2\}$, we get

\[
\frac{\omega_{A_i}(d_{S_n})}{n} \leq \omega_{A_i}(d) \leq \frac{\omega_{A_i}(d_{S_n})}{n-2-\alpha}.
\]

Recall that, by assumption, we have
\[
\omega_{A_1}(d)-c \cdot \omega_{A_2}(d) >0.
\]
One the other hand, we have, for every $n$,
\[
\frac{\omega_{A_1}(d_{S_n})}{n-2-\alpha}-c \cdot\frac{\omega_{A_2}(d_{S_n})}{n-2-\alpha} \geq \omega_{A_1}(d) - c \cdot\frac{\omega_{A_2}(d_{S_n})}{n} \cdot \frac{n}{n-2-\alpha} \geq \omega_{A_1}(d) - c \cdot\omega_{A_2}(d) \frac{n}{n-2-\alpha}.
\]
The right-hand-side is larger than $0$ for large enough $n$.
Hence, $\omega_{A_1}(d_{S_n})-c \cdot\omega_{A_2}(d_{S_n})>0$.
\end{proof}

\begin{remark}
Note that one can also obtain similar bounds in the opposite direction
\[
\omega_{A_1}(d)-c \cdot\omega_{A_2}(d) \geq \frac{\omega_{A_1}(d_{S_n})}{n} - c \cdot\frac{\omega_{A_2}(d_{S_n})}{n-2-\alpha} = \frac{1}{n} \left(\omega_{A_1}(d_{S_n}) - c \cdot\frac{n}{n-2-\alpha}\omega_{A_2}(d_{S_n})\right).
\]
One could hope to deduce from these exponential genericity of the limit action $G \curvearrowright X$ if we knew that the quantity on the right hand side is positive for some $n$ the approximant metrics $d_{S_n}$.
\label{rmk:oppositebound}
\end{remark}
 We now give the proof of Proposition~\ref{prop:gapimpliesexpgenericity}.

\begin{proof}[Proof of Proposition~\ref{prop:gapimpliesexpgenericity}] Since $\omega_G(d)$ is a true limit, the fact that $\omega_G(d)> \omega_A(d)$ gives 

\[\lim_{k \rightarrow \infty} \frac{\log \#\{g \in G: d(x_0,g x_0) \leq k\}}{k}> \limsup_{k \rightarrow \infty} \frac{\log \#\{g \in G: g \notin A \text{ and }d(x_0,g x_0) \leq k\}}{k},\] and in particular, given $\epsilon>0,$ there exists $n_0$ such that for all $k \geq n_0$ we have

$$\frac{\log \#\{g \in G: d(x_0,g x_0) \leq k\}}{k}>\frac{\log \#\{g \in G: g \in A \text{ and } d(x_0,g x_0) \leq k\}}{k}+\epsilon.$$ Multiplying both sides by $k$ and taking exponentials gives us $$\#\{g \in G: d(x_0,g x_0) \leq k\} > \#\{g \in G: g \in A \text{ and } d(x_0,g x_0) \leq k\}.e^{\epsilon.k}.$$ This yields one of the implications. The other implication follows by starting from the last inequality, and going backwards (in particular, it does not use the true limit assumption).
\end{proof}

\subsection{Exponential genericity} From the above, the proof of Theorem~\ref{thm:genericity} and Theorem~\ref{thm:p.A are dense} follow readily.

\begin{proof}[{Proof of Theorem~\ref{thm:genericity}}]
Let $(X,d)$ be the geodesic metric space with strongly contracting elements on which $G$ acts, and let $A$ be the collection of elements which are not strongly contracting on $(X,d)$. We have $\omega_A(d)<\omega_G(d)$ (by ~\cite{Yang2018}). Hence, using Theorem \ref{thm:Most_main_theorem}, we obtain that $\omega_{A}(d_{S_n})<\omega_{G}(d_{S_n})$ for all $n$ large enough. By Proposition~\ref{prop:rate_to_density}, since $\omega_G(d_{S_n})$ is a true limit, we have exponential genericity.
\end{proof}

\begin{proof}[{Proof of Theorem~\ref{thm:p.A are dense}}]
Let $(X,d)$ be the injective space on which $G=\MCG(\Sigma)$ acts geometrically.
In~\cite{SistoZalloum2024} it is shown that pseudo-Anosov mapping classes act as strongly contracting elements on these spaces (in fact, they are characterized by this property).
Hence, the result follows by Theorem~\ref{thm:genericity}.
\end{proof}

\subsection{Improving the hierarchies}

Now, we provide a proof to Theorem \ref{thm:good_hierarchies}.

\begin{proof} In \cite[Theorem~H]{PetytZalloum24}, Petyt and the third author construct for each mapping class group $\MCG(\Sigma)$ a geometric action on an injective space $(X,d)$ so that each $g.x,h.x$ are connected by a hierarhcy path which is also an $\alpha$-rough geodesic; where $\alpha$ depends only on $\Sigma$. Now, by Lemma \ref{lem:density}, we have $(n-1-\alpha)d_{S_n}(g,h)-(n-1) \leq d(gx_0,hx_0) \leq n d_{S_n}(g,h)$ for each large enough $n.$ Dividing by $n$ provides us with the following equation

$$\left(1-\frac{1+\alpha}{n}\right)d_{S_n} - \left(1-\frac{1}{n}\right) \leq \frac{d}{n} \leq d_{S_n}.$$

The metric $d_n:=\frac{d}{n}$ is injective for each $n$ and each $\alpha$-rough geodesic in $d$, becomes an $\frac{\alpha}{n}$-rough geodesic in $d_n$. The statement then easily follows by pushing the $\frac{\alpha}{n}$-rough geodesics to $d_{S_n}$ via the above equation.

\end{proof}

\subsection{Cogrowth} Finally we prove Theorem~\ref{thm:cogrowth_mcg}.

\begin{proof}[Proof of~Theorem~\ref{thm:cogrowth_mcg}]
Since $(X,d)$ admits an action with a strongly contracting element, it has pure exponential growth by~\cite{Yang17:ConvexCocompact}.
Hence, by \cite[Theorem~1.13]{CGZ24:Confined}, $\lambda_H(d) > \lambda_{G}(d)/2$. Hence, by Theorem~\ref{thm:Most_main_theorem}, the result follows.
\end{proof}

\section{Proof of growth tightness}
\label{sec:proof_GrowthTightness}

In this section we give the proof of the growth tightness result, Theorem~\ref{mainthm:cogrowth}, which we reestate here again for convenience.

\begin{restatethm}{GrowthTight4SomeS}
Let $G$ act properly coboundedly on a geodesic space $(X,d)$ with a strongly contracting element, and let $x_0 \in X$. Then, there exists $n_0 > 0$ such that for all $n \geq n_0$, if $S_n$ consists of elements whose $x_0$-orbit lies in $B(x_0,n) \subset X$, then $S_n$ generates $G$ and satisfies 
\[
\omega_{G/H}(d_{\bar{S}}) < \omega_G(d_S)
\]
for any infinite confined subgroup $H \subset G$ with a non-degenerate confining subset.
\end{restatethm}

This argument is more involved that the line of ideas in Section~\ref{sec:proofs}, since it entails checking infinitely many strict inequalities of the type $\omega_A(d_{S}) < \omega_G(d_{S})$, one for each infinite index normal subgroup $H$ of $G$. Unfortunately, one cannot hope to get a uniform  (in $H$) lower bound on the gap, since even in the free group case there are examples of ~\cite[Corollary~1.4]{JM20:FreeGap} sequences of subgroups $H_n$ with gap tending to $0$ (in the standard word-metric of $F_n$).

In order to prove growth tightness, we need a finer control on how the approximating word-metrics $d_{S_n}$ in the Cayley graph and the strongly contracting metric compare.

\subsection{Sketch of the Proof}
\label{subsec:sketch}
To enhance clarity, we begin with an overall sketch of the proof.

 The primary goal is to construct a generating set $S$ with the desired property that is independent of the confined subgroup $H$. Given an element $(X, d) \in \DDG$ with contracting elements, we obtain $S=S_n$ by selecting a large ball of radius $n$ in $X$, ensuring that the quasi-geodesic constants $(K_n, C_n)$ for the orbit map from $\Cay(G, S)$ to $X$ satisfy that $K_n$ is sufficiently close to $1$. The required closedness to $1$ depends on $\omega_X(G)$ and is mediated by the constant $a$ in Lemma~\ref{projLip} in Section~\ref{subsec:projection}. 
We furthermore fix a finite set of independent elements $\mathcal{F}$, chosen once and for all, which does not depend on the choice of $S$ or $H$, but only depends on the metric $d$.  Our goal will be to show that $\omega_{G/H}(d_{S_n}) < \omega_{G}(d_{S_n})$ for all $n$ large enough.

A key tool in our proof will be the following \emph{gap criterion} (Section~\ref{subsec:gap}).
Roughly, for a subset $A \subset G$, it gives conditions to construct, for each element $g\in A$, and constant $L>0$, a set of $2^{\| g \|_S/L}$ many \emph{different} elements in $G$ of length approximately $\|g\|_S$. A computation using critical exponent series then yields $\omega_A(S) < \omega_G(S)$.
More precisely, let $\mathcal{P}(n,L)$ denote the subsets of integer multiples of $L$ ranging from $L$ to $\lfloor n/L \rfloor L$, where multiples of $L^2$ are enforced (See Definition \ref{defPSet}). 
The gap criterion says that if for a given $A \subset G$,
we can construct:
\begin{enumerate} \item for every $g \in A$, an injective map $\Phi_g \colon \mathcal{P}(|g|,L) \to G$ \item so that $g_1 \neq g_2$ implies $\mathrm{Im}(\Phi_{g_1}) \cap \mathrm{Im}(\Phi_{g_2}) = \emptyset$, 
\item and so that the word length (wrt to $S$) of $\Phi_g(I)$ is comparable to that of $g$ (up to a logarithmic error in $L$ and factors depending on $|I|$ and uniform constants fixed in the beginning and depending on $\omega_G(S)$),
\end{enumerate}
then it follows that $\omega_A(S) < \omega_G(S)$.

In practice, we will construct a map $\Phi_g$ satisfying condition (1) in the gap criterion, for every arbitrary element $g \in G$, and such that $\mathrm{Im}(\Phi_g) \subset Hg$ (see~ Section~\ref{subsec:extensionmap}). 

Specifically, $\Phi_g$ is constructed by inserting contracting elements into a word representing $g$ (in the word metric $d_S$) at positions specified by $I \in \mathcal{P}(\|g\|_S, L)$. We will call this $\Phi_g$ the \emph{extension map}, since it is construct by ``extending'' $g$ in multiple ways.

The injectivity of $\Phi_g$ is established in Section~\ref{subsec:injectivity}. This requires the $C$-contracting axes of elements in $\mathcal{F}$ to be sufficiently long, specifically exceeding $2a \log L + D$. Here $a$ is the constant alluded above in the construction of the generating set, appearing in Lemma~\ref{projLip}. Hence, it depends only on $S$ and on $d$, and is in particular independent of $H$. On the other hand, $D$ does depend on $H$. 

    The injectivity of the map follows from the fact that the inserted contracting axes are sufficiently long to persist after insertion. Specifically, we can detect where we inserted by examining the fellow-traveling contracting axes of the geodesic from $1$ to $\Phi_g(I)$.
    
    To elaborate, for each axis we insert, the fellow-travel property of admissible path ensures that the projections of $1$ and $\Phi_g(I)$ onto the axis are widely separated. Conversely, for each non-inserted axis, Lemma \ref{projLip} and Lemma \ref{projToAxis} guarantee that the projections of $1$ and $\Phi_g(I)$ onto the axis are close, giving that the length of the axis exceeds $2a\log L + D$. Together, these two facts establish injectivity.

 Finally, we will apply the gap criterion when $A$ is a set of suitable representatives of $G/H$ (Section~\ref{subsec:extensionmap}). In this case, the fact that $\mathrm{Im}(\Phi_g) \subset Hg$ will ensure condition (2) in the gap criterion.

    The desired bound on the length of $\Phi_g(I)$ required by condition (3) follows then by the choice of $a$ and of $S$. Hence, the gap criterion goes through, and we obtain $\omega_{G/H}(d_S) < \omega_G(S)$, as desired, finishing the proof of Theorem~\ref{mainthm:cogrowth}.

\subsection{Preliminaries}\label{SecPrelim}

We first introduce and recall some basic concepts and results on diameters of projections to $C$-contracting subsets.

\subsubsection{Contracting geodesics} 
Let $(\U, d)$ be a geodesic metric space. Let $Z$ be a subset of $\U$ and $x$ be a point in $\U$.  By $d(x, Z)$ we mean the set-distance 
between $x$ and $Z$, \emph{i.e.} 
\[
d(x, Z) : = \inf \big \{ d(x, y): y \in Z \big \}. 
\]
Let $\delta>0$ be a large universal constant and 
\[ \pi_{Z}(x) \coloneqq \big \{ y\in Z: d(x, y) \le d(x, Z)+\delta \big \} \]
be the set of closest point projections from $x$ to $Z$. 
Then $\pi_{Z}(x)$ is non-empty. We refer to $\pi_{Z}(x) $ as the \emph{projection set} of $x$ to $Z$. Define $d_Z^{\pi}(x,y) \coloneqq \diam{\pi_Z(x) \cup \pi_Z(y)}$ for subset $Z\subset X$ and $x, y\in X$, where $\diam \cdot$ denotes the \emph{diameter} of a subset in the metric $d$.

\begin{defn} \label{Def:Contracting}
We say a subset $Z \subseteq X$ is \emph{$C$-contracting} for a constant $C>0$ if,
for all pairs of points $x, y \in X$, we have
\[
d(x, y) \leq d(x, Z) \quad  \Longrightarrow  \quad d_Z^{\pi}(x,y) \leq  C.
\]
Any such $C$ is called a \emph{contracting constant} for $Z$. A collection of $C$-contracting subsets shall be referred to as a \emph{$C$-contracting system}.
We note that $d^{\pi}_Y$ satisfies the triangular inequality: $$d^{\pi}_Y(x, z) \le d^{\pi}_Y(x, y) + d^{\pi}_Y(y, z)$$ for any $x, y, z\in X$.

An element $h\in \isom(\U)$ is called \textit{contracting} if it acts by translation on a contracting bi-infinite quasi-geodesic. Equivalently, the map $n\in \mathbb Z\longmapsto h^no$ is a quasi-geodesic with a contracting image.  
\end{defn}

Given a contracting element $h$, there exist a maximal elementary group $E(h)$ containing $\langle h \rangle$ as a finite index subgroup. Precisely,
\[
E(h) \coloneqq \{ g \in G : \exists n >0, (gh^ng^{-1}=h^n) \mbox{ or } (gh^ng^{-1}=h^{-n})  \}.
\]
We will define
\[
\ax(h) \coloneqq \{ f \cdot o | f \in E(h) \}
\]
and denote it the \emph{axis} of $h$.

\begin{defn}
	Let $\gamma: [0, 1]\to X$ be a path in $X$. Then $\gamma$ is called {\it $(K, C)$-quasi-geodesic} if for any $0\le a < b \le 1$, $length(\gamma|_{[a,b]})\le K d(\gamma(a), \gamma(b)) + C$ where $\gamma|_{[a,b]}$ is the subpath of $\gamma$ with endpoints $\gamma(a)$ and $\gamma(b)$.
\end{defn}
Let $G$ be a group with a finite generating set $S$ and Cayley graph $\Gamma$. Let $[g, h]_S$ be the geodesic in $\Gamma$ with endpoints $g$ and $h$. Denote the word metric on $\Gamma$ by $d_S$.

Fix a basepoint $o\in X$. Let $\pi: G\to X$ be defined by $g\mapsto go$. Then $\pi([g, h]_S)$ is a $(K, C)$-quasi-geodesic for some $K\ge 1$, $C\ge 0$.
Contracting subsets are Morse, and are preserved up to taking finite Hausdorff distance. The following properties shall be used later on.
\begin{lem}\label{BigFive}
	Let $Y\subseteq \U$ be a  $C$-contracting subset for $C>0$. Then the following holds.
	
	\begin{enumerate}
		\item
		For any  geodesic $\gamma$, we have $$\big |d^{\pi}_Y(\gamma^-,\gamma^+)- \diam{\gamma\cap N_C(Y)}\big|\leq 4C.$$
		
		\item
		There exists $\hat C=\hat C(C)$ such that   $d^{\pi}_Y(y, z)\le d(y, z)+\hat C$ for any $y,z\in \U$.
		
		\item 
		There exists $\hat C=\hat C(C)$ such that for any geodesic $\gamma$, we have $$\left| \diam{\pi_Y(\gamma)} - d^{\pi}_Y(\gamma_-, \gamma_+) \right| \leq \hat C$$
		
	\end{enumerate}
\end{lem}
We refer the reader to \cite{Yang:PS&growth} for a proof.

A $C$-contracting system $\mathcal{X}$ has the \emph{bounded intersection property} if there exist a constant $B \geq 0$ so that $\diam \pi_{X}(X') \leq B$ for every $X,X' \in \mathcal{X}$ with $X' \neq X$.
For example, it is easy to see that the collection $\{ g \ax(h) : g \in G\}$ is a contracting system with bounded intersection.
\begin{defn}
We say that a collection of contracting elements $\mathcal{F} = \{ f_i \}_{i \in \mathcal{I}} \subset G$ are \emph{independent} if 
for every $i, j \in \mathcal{I}$ with $i \neq j$, the collection 
\[
\{ g\ax(f_k) : g \in G, k\in \{i,j \}\}
\]
is a contracting system with bounded intersection.
\end{defn}

\subsubsection{Extension Lemma}
    In this paper, we frequently construct a path labeled by a word $(g_1, g_2,\cdots,g_n)$, which by convention means the following concatenation
$$
[o,g_1o]\cdot g_1[o,g_2o]\cdots (g_1\cdots g_{n-1})[o,g_no]
$$
where the basepoint $o$ is understood in context. With this convention, the paths labeled by $(g_1,g_2,g_3)$ and $(g_1g_2, g_3)$ respectively differ, depending on whether  $[o,g_1o]g_1[o,g_2o]$ is a geodesic or not.

{We fix a finite set $F \subseteq \Gamma$ of independent contracting elements and let $\f = \{g \ax(f) : f \in F, g \in \Gamma\}$.} The following notion of an admissible path allows     to construct   a quasi-geodesic  by concatenating geodesics via $\f$.

\begin{defn}[Admissible Path]\label{AdmDef} Given $L,\tau\geq0$, a path $\gamma$ is called $(L,\tau)$-\textit{admissible} in $\U$, if $\gamma$ is a concatenation of geodesics $p_0q_1p_1\cdots q_np_n$ $(n\in\mathbb{N})$, where the two endpoints of each $p_i$ lie in some $X_i\in \f$, and   the following   \textit{Long Local} and \textit{Bounded Projection} properties hold:
\begin{enumerate}
\item[(LL)] Each $p_i$  for $1\le i< n$ has length bigger than $L$, and  $p_0,p_n$ could be trivial;
\item[(BP)] For each $X_i$, we have $X_i\ne X_{i+1}$ and $\max\{\diam{\pi_{X_i}(q_i)},\diam{\pi_{X_i}(q_{i+1})}\}\leq\tau$, where $q_0:=\gamma^-$ and $q_{n+1}:=\gamma^+$ by convention.
\end{enumerate} 
The collection of sets $\{X_i: 1\le i\le n\}$ will be referred to as the \emph{ contracting subsets} associated to the admissible path.
\end{defn}

\begin{remark}\label{ConcatenationAdmPath}
    The path $q_i$ could be allowed to be trivial, so by the (BP) condition, it suffices to check $X_i\ne X_{i+1}$. It will be useful to note that admissible paths could be concatenated as follows: Let $p_0q_1p_1\cdots q_np_n$ and $p_0'q_1'p_1'\cdots q_n'p_n'$ be $(L,\tau)$-admissible. If $p_n=p_0'$ has length bigger than $L$, then the concatenation $(p_0q_1p_1\cdots q_np_n)\cdot (q_1'p_1'\cdots q_n'p_n')$ has a natural $(L,\tau)$-admissible structure.  
\end{remark}

    We frequently use the following lemma to show $X_i\neq X_{i+1}$ in property (BP):
    \begin{lem}\label{LL2ImplyBP}
        Using notations in the definition \ref{AdmDef} of admissible path, for all $1\le i\le n$, if\\ $\max\{ \diam{\pi_{X_{i}}(q_{i})}, \diam{\pi_{X_{i-1}}(q_{i})}\}\leq\tau$ and $\ell(q_{i}) > \tau$, then $X_i \neq X_{i-1}$.
    \end{lem}
    \begin{proof}
        Suppose to the contrary that $X_i = X_{i-1}$. Let $a = (q_i)^-$ and $b = (q_i)^+$. Then $a = (p_{i-1})^+$, $b = (p_i)^-$. Since $p_i\subset X_i$, $p_{i-1}\subset X_{i-1}$, we have $a, b\in X_i$. But $d(a, b) = d(\Pi_{X_i}(a), \Pi_{X_i}(b)) < \diam{\Pi_{X_i}(q_i)} \le \tau$. Since $a,b$ are the endpoints of $q_i$, $\ell(q_i)\le \tau$. This contradicts our assumption that $\ell(q_i) > \tau$.
    \end{proof}

A sequence of points $x_i$ in a path $p$   is called \textit{linearly ordered} if $x_{i+1}\in [x_i, p^+]_p$ for each $i$.

\begin{defn}[Fellow travel]\label{Fellow}
Let   $\gamma = p_0 q_1 p_1 \cdots q_n p_n$ be an $(L, \tau)-$admissible
path. We say $\gamma$ has \textit{$r$-fellow travel} property for some $r>0$   if for any geodesic  
$\alpha$  with the same endpoints as $\gamma$,   there exists a sequence of linearly ordered points $z_i,
w_i$ ($0 \le i \le n$) on $\alpha$ such that  
$$d(z_i, p_{i}^-) \le r,\quad d(w_i, p_{i}^+) \le r.$$
In particular, $\diam {N_r(X_i)\cap \alpha} \ge L$ for each $X_i\in \f(\gamma)$. 
\end{defn}
The following result  says that   a local long admissible path enjoys the fellow travel property.

\begin{prop}\label{admisProp}\cite{Yan14:GrowthTight}
For any $\tau>0$, there exist $L,  r, c>0$ depending only on $\tau,C$ such that  any $(L, \tau)$-admissible path   has $r$-fellow travel property. In particular, it is a $c$-quasi-geodesic.
\end{prop}

The next lemma gives a way to build admissible paths.
\begin{lem}[Extension Lemma]\label{extend3}

{For any independent contracting elements $h_{1}, h_{2}, h_{3} \in G$,} there exist constants  $L, r, B>0$ depending only on $C$ with the following property.  

Choose any element $f_i\in \langle h_i\rangle$ for each $1\le i\le 3$  to form the set $F$ of $f_i$'s satisfying $\min \{\| f_i \|_X \} \geq L$. Let $g,h\in G$ be any two elements.

There exists an element $f \in F$ such that   the path  $$\gamma:=[o, go]\cdot(g[o, fo])\cdot(gf[o,ho])$$ is an $(L, \tau)$-admissible path relative to $\f$. 
\end{lem}

\begin{remark} 
Since admissible paths are local conditions, we can connect via $F$  any number of elements  $g\in G$ to satisfy (1) and (2). We refer the reader to \cite{Yan14:GrowthTight} for a precise formulation.
\end{remark}

\begin{cor}[{\cite[Corollary~3.9]{Yan14:GrowthTight}}]\label{ProjMidSet}
    Let $\gamma$ be an admissible path and $X_k$ be a contracting subset for $\gamma$ where $0\le k\le n$. Then there exists $N = N(C, \tau) > 0$ such that $\diam{\Pi_{X_k}([\gamma_{-}, (p_k)_{-}]_{\gamma})} < N$ where $\gamma_{-}$ and $(p_k)_{-}$ denote the starting endpoints of $\gamma$ and $p_k$ respectively.
\end{cor}
The following lemma allows us to find a contracting element that is independent of the non-degenerate confining subset.
\begin{lem}[\cite{DY24:GrowthTight}]\label{lem:indepsetF}
    Let $P$  be a finite non-degenerate set. Then there exists a finite set $F$ of independent contracting elements in $G$ so that for any $p\in P$ and $f\in F$, we have $p \ax(f)\ne \ax(f)$.      
\end{lem}
\subsection{Projection of quasi geodesics to contracting subsets}
\label{subsec:projection}
The results in this section will allow us to bound the diameter of the projection of quasi-geodesic to contracting subsets of $Y$. 
The first result shows the bound on the diameter can be improved arbitrarily as the multiplicative constant $K$ in the quasi-isometry goes to $1$.

\begin{lem}\label{projLip}
	For any $\beta>0$, $C\ge 0, K\ge 1$, there exist $a = a(K, \beta) \ge 0$, $b = b(\beta, K, C) > 0$ and $c = c(K) > 0$ with the following property. 
	Let $Y$ be a $\beta$-contracting subset of $X$ and $\gamma$ be a $(K, C)$-quasi-geodesic with endpoints $x\in X$ and $y\in Y$. Then $$\diam {\pi_{Y}(\gamma)} \le a \log d(x, y) + b + c d(\pi_Y(x), y)).$$ Furthermore,  we may choose $a(1, \beta)=0$ and $\lim_{K\to 1} a(K, \beta) = 0$ for any $\beta > 0$. 
\end{lem}
\begin{proof}
	Set $a(1, \beta)=0$ and assume that $K>1$. Denote $\theta \coloneqq 1-\frac{1}{K}$.
	Let $z$ be the first point on $\gamma$ such that $d(z, Y) \le  \theta d(x, Y)$ and let $\gamma_1$ be the subpath of $\gamma$ between $x$ and $z$. Let $\mathfrak h = \max\left\{\frac{2K^2\beta}{K-1}, \frac{CK}{K^2-1}, 1\right\}$.
	\begin{claim}
		If $d(x, Y) > \mathfrak h$, then $\diam{\pi_Y(\gamma_1)} \le 6(K + 1) \beta$.
	\end{claim}
	\begin{proof}[Proof of the claim]
		Let $D = \diam{\pi_Y(\gamma_1)}$, $l = \mathrm{length}(\gamma_1)$ and $h = d(x, Y)$. Then by triangular inequality, $$d(x, z) \le d(x, Y) + d(z, Y) + d^{\pi}_Y(x, z) \le h + \theta h + D.$$ 
		Since $\gamma$ is a $(K, C)$-quasi-geodesic, we have \begin{align}\label{lboundEq}
			l \le K d(x, z) + C \le K(h(1+\theta) + D) + C
		\end{align}
		
		We are going to bound $D$. To that end,  we write $\gamma_1 = \gamma_2 \cup \gamma_3$ as the concatenation of two paths where  the subpath $\gamma_2$ has length $h$ starting at $x$ and thus is contained in $B(x, h)$, and $\gamma_3$ is the remaining subpath of $\gamma_1$ with $\mathrm{length}(\gamma_3) = l - h$. Since $Y$ is $\beta$-contracting and $h = d(x, Y)$, we have $\diam{\pi_Y(\gamma_2)} \le \beta$. For $\gamma_3$, we may divide it into disjoint subpaths of length $\theta h$, and the number of them is at most $({l-h})/{\theta h} + 1.$ Since the path $\gamma_1$ from $x$ to $z$ has distance at least $\theta h$ to $Y$, we have that each of these subpaths is contained in a ball which is disjoint from $Y$. Thus, $$\diam{\pi_Y(\gamma_3)}\le \beta\left({(l-h)}/{\theta h} + 1\right).$$ Taking into account the projections  of $\gamma_2$ and $\gamma_3$, we have $$D = \diam{\pi_Y(\gamma_1)} \le \diam{\pi_Y(\gamma_2)} + \diam{\pi_Y(\gamma_3)} \le \beta \left(2 + {(l-h)}/{\theta h}\right)$$
		which yields $$l \ge h + \left({D}/{\beta} - 2\right)\theta h.$$
		With (\ref{lboundEq}) we have $$h + \left({D}/{\beta} - 2\right) \theta h \le K\left(h(1+\theta) + D\right) + C,$$ or equivalently, using  $\theta=1-\frac{1}{K}$, $$D\left( {\theta h}/{\beta} - K\right)  \le h\left(K +K\theta- 1+2\theta \right) + C =  2h(K +1 )\theta + C.$$ Recall that   $\mathfrak h= \max\{\frac{2K\beta}{\theta}, \frac{C}{(K+1)\theta}, 1\}$. If $h > \mathfrak h$, we have ${h\theta }/{\beta} \ge 2K$ and $C\le h (K +1) \theta$. Thus, $$D \le \frac{3h(K +1)\theta}{{\theta h}/{2\beta}} = 6(K + 1)\beta.$$
		This proves the claim.
	\end{proof}
	Let $x_0 = x$ and define inductively $x_k$ to be the first point on $\gamma$ with $d(x_k, Y) \le \theta d(x_{k-1}, Y)$ for $k\ge 1$. Let $m$ be the minimal $k$ such that $d(x_k, Y)\le \mathfrak h$. Then $$-{\log_{\theta} \frac{d(x,Y)}{\mathfrak h}} \le m \le 1- {\log_\theta \frac{d(x,Y)}{\mathfrak h}} \le 1 - \log_\theta d(x,Y).$$  
	Let $\gamma_k$ be the subpath of $\gamma$ between $x_{k-1}$ and $x_k$. As    $d(\gamma_{k}, Y)\ge \mathfrak h$ for $k\le m$, by the claim above, we have $\diam{\pi_Y(\gamma_k)}\le 6(K+1) \beta$. If $\alpha_1$ denotes the subpath of $\gamma$ between $x$ and $x_m$, then $$\diam{\pi_Y(\alpha_1)} \le 6m(K+1) C,$$ and especially $d_Y^\pi(x, x_m) \le 6m(K+1) \beta$. Thus, $d(x_m, y)\le d(x_m, \pi_Y(x_m)) + d_Y^\pi(x, x_m) + d(\pi_Y(x), y) \le 6m(K+1)\beta + \mathfrak h + d(\pi_Y(x), y)$. Let $\alpha_2$ be the subpath of $\gamma$ between $x_m$ and $y$. Since $\alpha_2$ is also a $(K, C)$-quasi-geodesic, we have $\mathrm{length}(\alpha_2)\le K d(x_m, y) + C$. So by Lemma \ref{BigFive}, $$\diam{\pi_Y(\alpha_2)} \le \mathrm{length}(\alpha_2) + \hat C \le 6m K(K+1) \beta + K \mathfrak h + K d(\pi_Y(x), y) + C + \hat C.$$ Finally, $$\begin{aligned}
		\diam{\pi_Y(\gamma)} &\le \diam{\pi_Y(\alpha_1)} + \diam{\pi_Y(\alpha_2)}\\
		&\le 6m (K+1)^2 \beta + K \mathfrak h + K d(\pi_Y(x), y) + C + \hat C\\
		&\le a \log d(x, y)  +b + c d(\pi_Y(x), y) 
	\end{aligned}$$
	where $a = -{6(K+1)^2 \beta}/{\log \theta}$, $b = {6(K+1)^2 \beta} + K \mathfrak h + C + \hat C$ and $c = K$. 

    Noticing  $\lim_{K\to 1^+} a = 0$, the lemma thus follows.
\end{proof}

\subsection{Construction of Extension Map}
\label{subsec:extensionmap}
In this section we will construct, for any $g \in G$ the ``extension map'' $\Phi_g \colon \mathcal{P}(n, L) \to Hg$ required in the ``gap criterion''.
This section follows~\cite{DY24:GrowthTight} closely. The main difference lies in that here the choice of the contracting axis depends on the local structure of $I$.

\begin{lem}\label{BddProjLem}
	Let $\mathcal{F}_0$ be an infinite set of   independent contracting elements in $G$. For any integer $n\ge 1$, there exist a finite subset $\mathcal{F}$ in  $\mathcal{F}_0$   and a constant $B>0$ with the following property. If $(a_1,\cdots, a_n)$ is any $n$-tuple in $G$, then there exists $f\in \mathcal F$ such that $\diam{\pi_{\ax (f)}([o, a_i o])} \le B$ for each $1\le i\le n$.
\end{lem}

This lemma follows directly from~\cite[Lemma 2.12 and Lemma 2.17]{Yang17:ConvexCocompact}.

In Lemma~\ref{BddProjLem}, we pick $n=5$ and get a finite set $\mathcal{F}$ consisting of contracting element $f$, and the choice of $\mathcal{F}$ is independent of the generating set $S$ and of the normal subgroup $H$. This set will be fixed once and for all.

Now, let $P$ be a finite confining subset for the confined subgroup $H$.
Assume that $P$ is non-degenerate, i.e. disjoint with $E(G)$. By lemma~\ref{lem:indepsetF}, we pick a contracting element $f\in G$ such that $\diam{\pi_{\ax(f)}([o, po])} \le B$ for any $p\in P$. Here $f$ depends on $P$.

Recall that by Proposition~\ref{admisProp}, there exist $\mathcal{D}>0, \lambda>1, r>0$ such that any $(\mathcal{D}, B)$-admissible path of $\max\{C, C'\}$-contracting subsets has \emph{$r$-fellow travel} property and is an $(\lambda, \lambda)$-quasi-geodesic. By taking a high power, we may assume that $\|f\|_X > \mathcal{D}$ and $\|f_i\|_X > \mathcal{D}$ for each $f_i\in \mathcal{F}$.

The following natural definition will be important in order to get enough control when approximating $X$ by generating sets.

\begin{defn}\label{defPSet}
	Let $L\in \mathbb{N}$. Define $\mathcal{P}(n, L)$ to be the  family of non-empty subsets $I\subset\{L, 2L \cdots, \lfloor \frac{n}{L}\rfloor L\}$ such that $iL^2\in I$ for each $1\le i\le \lfloor \frac{n}{L^2}\rfloor$. 
\end{defn}

Let $g\in G$ be an arbitrary element of $G$. If $\|g\|_S=n$, $s_i \in S$ for $i=1,\cdots,n$, and we have $g=s_1\cdots s_n$ which is a word geodesic of $g$. Denote $g_0 = id$ and $g_i = s_1 \cdots s_i$ for $1\le i\le n$.

This construction modifies the \emph{extension map} in \cite{DY24:GrowthTight}. Roughly, for each element $g \in G$, this map $\Phi_g$ yields $2^{\lfloor n/L \rfloor}$ many elements in the coset of $g$, each with length at most $n(1 + \frac{1}{L})$. It does so by ``extending'' $g$ inserting words from axes of contracting elements at prescribed places (see Figure~\ref{fig:Con_Ext} for an illustration).

\begin{defn}\label{defExtensionMap}
	Let $g \in G$, $g=   s_1 \cdots s_n\in G$ be a word geodesic and $L\ge 1$. The {\it extension map} $\Phi_g: \mathcal P(n,L)\to G$ is  defined by $$I=(a_1, \cdots, a_m)\mapsto \left(\prod_{i=1}^m g_{a_{i-1}}^{-1} g_{a_i} f_i f p_i f^{-1} f_i^{-1}\right) g_{a_{m}}^{-1} g_n.$$
	where $f_i\in \mathcal F$ are chosen with the following properties:
	\begin{enumerate}\label{contrachoice}
		\item\label{prevproj}  $[g_{a_{i-1}} o, g_{a_i} o]$ has $B$-bounded projection to $g_{a_i}\ax(f_i)$;
		\item\label{nextproj} $[g_{a_{i}} o, g_{a_{i+1}} o]$ has $B$-bounded projection to $g_{a_i}\ax(f_i)$;
		\item\label{prevLproj} $[g_{a_{i} - L^2} o, g_{a_{i}} o]$ has $B$-bounded projection to $g_{a_i}\ax(f_i)$;
		\item\label{nextLproj} $[g_{a_{i}} o, g_{L^2 + a_{i}} o]$ has $B$-bounded projection to $g_{a_i}\ax(f_i)$;
        \item\label{fproj} $[o, f o]$ has $B$-bounded projection to $\ax(f_i)$;
	\end{enumerate}
	where    $a_i = 0$ for $i\le 0$, $a_{j} = n$ for $j\ge m+1$ by convention, and
	\begin{enumerate}\setcounter{enumi}{5}\label{confinedchoice} 
		\item $p_i\in P$ are chosen in the confining subset $P$ so that $g_{a_i}f_i f p_i(g_{a_i}f_if)^{-1}$ is contained in $H$.
	\end{enumerate} 
	
	Let $\gamma_g(I)$ be the path labeled by $$(g_{a_1}, f_1,1, f, p_1, f^{-1}, 1, f_1^{-1}, \cdots, g_{a_{i-1}}^{-1} g_{a_i}, f_i, 1, f, p_i, f^{-1}, 1, f_i^{-1}, \cdots, g_{a_{m}}^{-1} g).$$
\end{defn}
\begin{figure}
    \centering
    \tikzset{every picture/.style={line width=0.75pt}} 

    \begin{tikzpicture}[x=0.75pt,y=0.75pt,yscale=-1,xscale=1]

    \draw    (134,91) -- (176,91) ;
    \draw [shift={(178,91)}, rotate = 180] [color={rgb, 255:red, 0; green, 0; blue, 0 }  ][line width=0.75]    (10.93,-3.29) .. controls (6.95,-1.4) and (3.31,-0.3) .. (0,0) .. controls (3.31,0.3) and (6.95,1.4) .. (10.93,3.29) ;
    \draw [shift={(134,91)}, rotate = 89.28] [color={rgb, 255:red, 0; green, 0; blue, 0 }  ,draw opacity=1 ][fill={rgb, 255:red, 0; green, 0; blue, 0 }  ,fill opacity=1 ][line width=1.5]      (0, 0) circle [x radius= 2.36, y radius= 2.36]   ;
    \draw    (178,91) -- (220,91) ;
    \draw [shift={(222,91)}, rotate = 180] [color={rgb, 255:red, 0; green, 0; blue, 0 }  ][line width=0.75]    (10.93,-3.29) .. controls (6.95,-1.4) and (3.31,-0.3) .. (0,0) .. controls (3.31,0.3) and (6.95,1.4) .. (10.93,3.29)   ;
    \draw    (222,91) -- (264,91) ;
    \draw [shift={(266,91)}, rotate = 180] [color={rgb, 255:red, 0; green, 0; blue, 0 }  ][line width=0.75]    (10.93,-3.29) .. controls (6.95,-1.4) and (3.31,-0.3) .. (0,0) .. controls (3.31,0.3) and (6.95,1.4) .. (10.93,3.29)   ;
    \draw    (310,91) -- (352,91) ;
    \draw [shift={(354,91)}, rotate = 180] [color={rgb, 255:red, 0; green, 0; blue, 0 }  ][line width=0.75]    (10.93,-3.29) .. controls (6.95,-1.4) and (3.31,-0.3) .. (0,0) .. controls (3.31,0.3) and (6.95,1.4) .. (10.93,3.29)   ;
    \draw [shift={(354,91)}, rotate = 89.28] [color={rgb, 255:red, 0; green, 0; blue, 0 }  ,draw opacity=1 ][fill={rgb, 255:red, 0; green, 0; blue, 0 }  ,fill opacity=1 ][line width=1.5]      (0, 0) circle [x radius= 2.36, y radius= 2.36]   ;
    \draw [color={rgb, 255:red, 0; green, 114; blue, 255 }  ,draw opacity=1 ]   (354,91) -- (396,91) ;
    \draw [shift={(398,91)}, rotate = 180] [color={rgb, 255:red, 0; green, 114; blue, 255 }  ,draw opacity=1 ][line width=0.75]    (10.93,-3.29) .. controls (6.95,-1.4) and (3.31,-0.3) .. (0,0) .. controls (3.31,0.3) and (6.95,1.4) .. (10.93,3.29)   ;
    \draw [color={rgb, 255:red, 255; green, 3; blue, 32 }  ,draw opacity=1 ]   (442,91) -- (441.06,123) ;
    \draw [shift={(441,125)}, rotate = 271.68] [color={rgb, 255:red, 255; green, 3; blue, 32 }  ,draw opacity=1 ][line width=0.75]    (10.93,-3.29) .. controls (6.95,-1.4) and (3.31,-0.3) .. (0,0) .. controls (3.31,0.3) and (6.95,1.4) .. (10.93,3.29)   ;
    \draw [color={rgb, 255:red, 0; green, 114; blue, 255 }  ,draw opacity=1 ]   (397,125) -- (355,125) ;
    \draw [shift={(353,125)}, rotate = 360] [color={rgb, 255:red, 0; green, 114; blue, 255 }  ,draw opacity=1 ][line width=0.75]    (10.93,-3.29) .. controls (6.95,-1.4) and (3.31,-0.3) .. (0,0) .. controls (3.31,0.3) and (6.95,1.4) .. (10.93,3.29)   ;
    \draw    (133,125) -- (175,125) ;
    \draw [shift={(177,125)}, rotate = 180] [color={rgb, 255:red, 0; green, 0; blue, 0 }  ][line width=0.75]    (10.93,-3.29) .. controls (6.95,-1.4) and (3.31,-0.3) .. (0,0) .. controls (3.31,0.3) and (6.95,1.4) .. (10.93,3.29)   ;
    \draw [shift={(133,125)}, rotate = 89.28] [color={rgb, 255:red, 0; green, 0; blue, 0 }  ,draw opacity=1 ][fill={rgb, 255:red, 0; green, 0; blue, 0 }  ,fill opacity=1 ][line width=1.5]      (0, 0) circle [x radius= 2.36, y radius= 2.36]   ;
    \draw    (177,125) -- (219,125) ;
    \draw [shift={(221,125)}, rotate = 180] [color={rgb, 255:red, 0; green, 0; blue, 0 }  ][line width=0.75]    (10.93,-3.29) .. controls (6.95,-1.4) and (3.31,-0.3) .. (0,0) .. controls (3.31,0.3) and (6.95,1.4) .. (10.93,3.29)   ;
    \draw    (221,125) -- (263,125) ;
    \draw [shift={(265,125)}, rotate = 180] [color={rgb, 255:red, 0; green, 0; blue, 0 }  ][line width=0.75]    (10.93,-3.29) .. controls (6.95,-1.4) and (3.31,-0.3) .. (0,0) .. controls (3.31,0.3) and (6.95,1.4) .. (10.93,3.29)   ;
    \draw    (265,125) -- (307,125) ;
    \draw [shift={(309,125)}, rotate = 180] [color={rgb, 255:red, 0; green, 0; blue, 0 }  ][line width=0.75]    (10.93,-3.29) .. controls (6.95,-1.4) and (3.31,-0.3) .. (0,0) .. controls (3.31,0.3) and (6.95,1.4) .. (10.93,3.29)   ;
    \draw [color={rgb, 255:red, 0; green, 114; blue, 255 }  ,draw opacity=1 ]   (398,91) -- (440,91) ;
    \draw [shift={(442,91)}, rotate = 180] [color={rgb, 255:red, 0; green, 114; blue, 255 }  ,draw opacity=1 ][line width=0.75]    (10.93,-3.29) .. controls (6.95,-1.4) and (3.31,-0.3) .. (0,0) .. controls (3.31,0.3) and (6.95,1.4) .. (10.93,3.29)   ;
    \draw [color={rgb, 255:red, 0; green, 114; blue, 255 }  ,draw opacity=1 ]   (441,125) -- (399,125) ;
    \draw [shift={(397,125)}, rotate = 360] [color={rgb, 255:red, 0; green, 114; blue, 255 }  ,draw opacity=1 ][line width=0.75]    (10.93,-3.29) .. controls (6.95,-1.4) and (3.31,-0.3) .. (0,0) .. controls (3.31,0.3) and (6.95,1.4) .. (10.93,3.29)   ;
    \draw    (353,125) -- (353,191) ;
    \draw [shift={(353,193)}, rotate = 270] [color={rgb, 255:red, 0; green, 0; blue, 0 }  ][line width=0.75]    (10.93,-3.29) .. controls (6.95,-1.4) and (3.31,-0.3) .. (0,0) .. controls (3.31,0.3) and (6.95,1.4) .. (10.93,3.29)   ;
    \draw [shift={(353,193)}, rotate = 89.28] [color={rgb, 255:red, 0; green, 0; blue, 0 }  ,draw opacity=1 ][fill={rgb, 255:red, 0; green, 0; blue, 0 }  ,fill opacity=1 ][line width=1.5]      (0, 0) circle [x radius= 2.36, y radius= 2.36]   ;
    \draw   (106,235) -- (123.5,235) -- (123.5,211) -- (158.5,211) -- (158.5,235) -- (176,235) -- (141,251) -- cycle ;
    \draw [color={rgb, 255:red, 0; green, 114; blue, 255 }  ,draw opacity=1 ]   (93.01,285) -- (93.01,326.05) ;
    \draw [shift={(93.01,328.05)}, rotate = 270] [color={rgb, 255:red, 0; green, 114; blue, 255 }  ,draw opacity=1 ][line width=0.75]    (10.93,-3.29) .. controls (6.95,-1.4) and (3.31,-0.3) .. (0,0) .. controls (3.31,0.3) and (6.95,1.4) .. (10.93,3.29)   ;
    \draw    (44,285) -- (91.01,285) ;
    \draw [shift={(93.01,285)}, rotate = 180] [color={rgb, 255:red, 0; green, 0; blue, 0 }  ][line width=0.75]    (10.93,-3.29) .. controls (6.95,-1.4) and (3.31,-0.3) .. (0,0) .. controls (3.31,0.3) and (6.95,1.4) .. (10.93,3.29)   ;
    \draw [shift={(42,285)}, rotate = 89.28] [color={rgb, 255:red, 0; green, 0; blue, 0 }  ,draw opacity=1 ][fill={rgb, 255:red, 0; green, 0; blue, 0 }  ,fill opacity=1 ][line width=1.5]      (0, 0) circle [x radius= 2.36, y radius= 2.36]   ;
    \draw [color={rgb, 255:red, 255; green, 3; blue, 32 }  ,draw opacity=1 ]   (93.01,328.05) -- (144,328.05) ;
    \draw [shift={(146,328.05)}, rotate = 180] [color={rgb, 255:red, 255; green, 3; blue, 32 }  ,draw opacity=1 ][line width=0.75]    (10.93,-3.29) .. controls (6.95,-1.4) and (3.31,-0.3) .. (0,0) .. controls (3.31,0.3) and (6.95,1.4) .. (10.93,3.29)   ;
    \draw    (146,285) -- (193.01,285) ;
    \draw [shift={(195.01,285)}, rotate = 180] [color={rgb, 255:red, 0; green, 0; blue, 0 }  ][line width=0.75]    (10.93,-3.29) .. controls (6.95,-1.4) and (3.31,-0.3) .. (0,0) .. controls (3.31,0.3) and (6.95,1.4) .. (10.93,3.29)   ;
    \draw    (248,285) -- (295.01,285) ;
    \draw [shift={(297.01,285)}, rotate = 180] [color={rgb, 255:red, 0; green, 0; blue, 0 }  ][line width=0.75]    (10.93,-3.29) .. controls (6.95,-1.4) and (3.31,-0.3) .. (0,0) .. controls (3.31,0.3) and (6.95,1.4) .. (10.93,3.29)   ;
    \draw  [dash pattern={on 4.5pt off 4.5pt}]  (297.01,285) -- (366.07,285) ;
    \draw    (366.07,285) -- (413.07,285) ;
    \draw [shift={(415.07,285)}, rotate = 180] [color={rgb, 255:red, 0; green, 0; blue, 0 }  ][line width=0.75]    (10.93,-3.29) .. controls (6.95,-1.4) and (3.31,-0.3) .. (0,0) .. controls (3.31,0.3) and (6.95,1.4) .. (10.93,3.29)   ;
    \draw    (468.07,285) -- (515.07,285) ;
    \draw [shift={(517.07,285)}, rotate = 180] [color={rgb, 255:red, 0; green, 0; blue, 0 }  ][line width=0.75]    (10.93,-3.29) .. controls (6.95,-1.4) and (3.31,-0.3) .. (0,0) .. controls (3.31,0.3) and (6.95,1.4) .. (10.93,3.29)   ;
    \draw [shift={(517.07,285)}, rotate = 89.28] [color={rgb, 255:red, 0; green, 0; blue, 0 }  ,draw opacity=1 ][fill={rgb, 255:red, 0; green, 0; blue, 0 }  ,fill opacity=1 ][line width=1.5]      (0, 0) circle [x radius= 2.36, y radius= 2.36]   ;
    \draw    (266,91) -- (308,91) ;
    \draw [shift={(310,91)}, rotate = 180] [color={rgb, 255:red, 0; green, 0; blue, 0 }  ][line width=0.75]    (10.93,-3.29) .. controls (6.95,-1.4) and (3.31,-0.3) .. (0,0) .. controls (3.31,0.3) and (6.95,1.4) .. (10.93,3.29)   ;
    \draw    (309,125) -- (351,125) ;
    \draw [shift={(353,125)}, rotate = 180] [color={rgb, 255:red, 0; green, 0; blue, 0 }  ][line width=0.75]    (10.93,-3.29) .. controls (6.95,-1.4) and (3.31,-0.3) .. (0,0) .. controls (3.31,0.3) and (6.95,1.4) .. (10.93,3.29)   ;
    \draw [shift={(353,125)}, rotate = 89.28] [color={rgb, 255:red, 0; green, 0; blue, 0 }  ,draw opacity=1 ][fill={rgb, 255:red, 0; green, 0; blue, 0 }  ,fill opacity=1 ][line width=1.5]      (0, 0) circle [x radius= 2.36, y radius= 2.36]   ;
    \draw [color={rgb, 255:red, 0; green, 114; blue, 255 }  ,draw opacity=1 ]   (146,285) -- (146,326.05) ;
    \draw [shift={(146,328.05)}, rotate = 270] [color={rgb, 255:red, 0; green, 114; blue, 255 }  ,draw opacity=1 ][line width=0.75]    (10.93,-3.29) .. controls (6.95,-1.4) and (3.31,-0.3) .. (0,0) .. controls (3.31,0.3) and (6.95,1.4) .. (10.93,3.29)   ;
    \draw [color={rgb, 255:red, 0; green, 114; blue, 255 }  ,draw opacity=1 ]   (195.01,285) -- (195.01,326.05) ;
    \draw [shift={(195.01,328.05)}, rotate = 270] [color={rgb, 255:red, 0; green, 114; blue, 255 }  ,draw opacity=1 ][line width=0.75]    (10.93,-3.29) .. controls (6.95,-1.4) and (3.31,-0.3) .. (0,0) .. controls (3.31,0.3) and (6.95,1.4) .. (10.93,3.29)   ;
    \draw [color={rgb, 255:red, 255; green, 3; blue, 32 }  ,draw opacity=1 ]   (195.01,328.05) -- (246,328.05) ;
    \draw [shift={(248,328.05)}, rotate = 180] [color={rgb, 255:red, 255; green, 3; blue, 32 }  ,draw opacity=1 ][line width=0.75]    (10.93,-3.29) .. controls (6.95,-1.4) and (3.31,-0.3) .. (0,0) .. controls (3.31,0.3) and (6.95,1.4) .. (10.93,3.29)   ;
    \draw [color={rgb, 255:red, 0; green, 114; blue, 255 }  ,draw opacity=1 ]   (248,285) -- (248,326.05) ;
    \draw [shift={(248,328.05)}, rotate = 270] [color={rgb, 255:red, 0; green, 114; blue, 255 }  ,draw opacity=1 ][line width=0.75]    (10.93,-3.29) .. controls (6.95,-1.4) and (3.31,-0.3) .. (0,0) .. controls (3.31,0.3) and (6.95,1.4) .. (10.93,3.29)   ;
    \draw [color={rgb, 255:red, 0; green, 114; blue, 255 }  ,draw opacity=1 ]   (415.07,285) -- (415.07,326.05) ;
    \draw [shift={(415.07,328.05)}, rotate = 270] [color={rgb, 255:red, 0; green, 114; blue, 255 }  ,draw opacity=1 ][line width=0.75]    (10.93,-3.29) .. controls (6.95,-1.4) and (3.31,-0.3) .. (0,0) .. controls (3.31,0.3) and (6.95,1.4) .. (10.93,3.29)   ;
    \draw [color={rgb, 255:red, 255; green, 3; blue, 32 }  ,draw opacity=1 ]   (415.07,328.05) -- (466.07,328.05) ;
    \draw [shift={(468.07,328.05)}, rotate = 180] [color={rgb, 255:red, 255; green, 3; blue, 32 }  ,draw opacity=1 ][line width=0.75]    (10.93,-3.29) .. controls (6.95,-1.4) and (3.31,-0.3) .. (0,0) .. controls (3.31,0.3) and (6.95,1.4) .. (10.93,3.29)   ;
    \draw [color={rgb, 255:red, 0; green, 114; blue, 255 }  ,draw opacity=1 ]   (468.07,285) -- (468.07,326.05) ;
    \draw [shift={(468.07,328.05)}, rotate = 270] [color={rgb, 255:red, 0; green, 114; blue, 255 }  ,draw opacity=1 ][line width=0.75]    (10.93,-3.29) .. controls (6.95,-1.4) and (3.31,-0.3) .. (0,0) .. controls (3.31,0.3) and (6.95,1.4) .. (10.93,3.29)   ;
    \draw    (354,25) -- (354,91) ;
    \draw [shift={(354,23)}, rotate = 90] [color={rgb, 255:red, 0; green, 0; blue, 0 }  ][line width=0.75]    (10.93,-3.29) .. controls (6.95,-1.4) and (3.31,-0.3) .. (0,0) .. controls (3.31,0.3) and (6.95,1.4) .. (10.93,3.29)   ;
    \draw [shift={(354,23)}, rotate = 89.28] [color={rgb, 255:red, 0; green, 0; blue, 0 }  ,draw opacity=1 ][fill={rgb, 255:red, 0; green, 0; blue, 0 }  ,fill opacity=1 ][line width=1.5]      (0, 0) circle [x radius= 2.36, y radius= 2.36]   ;
    \draw    (57,42) -- (352.03,90.67) ;
    \draw [shift={(354,91)}, rotate = 189.37] [color={rgb, 255:red, 0; green, 0; blue, 0 }  ][line width=0.75]    (10.93,-3.29) .. controls (6.95,-1.4) and (3.31,-0.3) .. (0,0) .. controls (3.31,0.3) and (6.95,1.4) .. (10.93,3.29)   ;
    \draw [shift={(55.03,41.67)}, rotate = 89.28] [color={rgb, 255:red, 0; green, 0; blue, 0 }  ,draw opacity=1 ][fill={rgb, 255:red, 0; green, 0; blue, 0 }  ,fill opacity=1 ][line width=1.5]      (0, 0) circle [x radius= 2.36, y radius= 2.36]   ;
    \draw    (353,125) -- (162.89,192.33) ;
    \draw [shift={(161,193)}, rotate = 340.5] [color={rgb, 255:red, 0; green, 0; blue, 0 }  ][line width=0.75]    (10.93,-3.29) .. controls (6.95,-1.4) and (3.31,-0.3) .. (0,0) .. controls (3.31,0.3) and (6.95,1.4) .. (10.93,3.29)   ;
    \draw [shift={(161,193)}, rotate = 89.28] [color={rgb, 255:red, 0; green, 0; blue, 0 }  ,draw opacity=1 ][fill={rgb, 255:red, 0; green, 0; blue, 0 }  ,fill opacity=1 ][line width=1.5]      (0, 0) circle [x radius= 2.36, y radius= 2.36]   ;

    \draw (98,76) node [anchor=north west][inner sep=0.75pt]   [align=left] {$\displaystyle g_{a_{i-1}} o$};
    \draw (324,70) node [anchor=north west][inner sep=0.75pt]   [align=left] {$\displaystyle g_{a_{i}} o$};
    \draw (177,76) node [anchor=north west][inner sep=0.75pt]   [align=left] {$\displaystyle s_{a_{i-1} +1} \cdots s_{a_{i}}$};
    \draw (363,93) node [anchor=north west][inner sep=0.75pt]   [align=left] {$\displaystyle f_{i}$};
    \draw (416,127) node [anchor=north west][inner sep=0.75pt]   [align=left] {$\displaystyle f$};
    \draw (377,128) node [anchor=north west][inner sep=0.75pt]   [align=left] {$\displaystyle f_{i}$};
    \draw (356,152) node [anchor=north west][inner sep=0.75pt]   [align=left] {$\displaystyle g_{a_{i}}^{-1} g_{a_{i+1}}$};
    \draw (356,48) node [anchor=north west][inner sep=0.75pt]   [align=left] {$\displaystyle g_{a_{i-1}}^{-1} g_{a_{i}}$};
    \draw (162,42) node [anchor=north west][inner sep=0.75pt]   [align=left] {$\displaystyle g_{a_{i}-L^2}^{-1} g_{a_{i}}$};
    \draw (405,93) node [anchor=north west][inner sep=0.75pt]   [align=left] {$\displaystyle f$};
    \draw (444,94) node [anchor=north west][inner sep=0.75pt]   [align=left] {$\displaystyle p_{i}$};
    \draw (51.99,263) node [anchor=north west][inner sep=0.75pt]   [align=left] {$\displaystyle g_{a_{1}}$};
    \draw (147.93,260) node [anchor=north west][inner sep=0.75pt]   [align=left] {$\displaystyle g_{a_{1}}^{-1} g_{a_{2}}$};
    \draw (105.01,286) node [anchor=north west][inner sep=0.75pt]   [align=left] {$\displaystyle f_{1} f$};
    \draw (108.26,330.2) node [anchor=north west][inner sep=0.75pt]   [align=left] {$\displaystyle p_{1}$};
    \draw (208.01,287) node [anchor=north west][inner sep=0.75pt]   [align=left] {$\displaystyle f_{2} f$};
    \draw (210.26,330.2) node [anchor=north west][inner sep=0.75pt]   [align=left] {$\displaystyle p_{2}$};
    \draw (427.01,287) node [anchor=north west][inner sep=0.75pt]   [align=left] {$\displaystyle f_{m} f$};
    \draw (430.26,330.2) node [anchor=north west][inner sep=0.75pt]   [align=left] {$\displaystyle p_{m}$};
    \draw (247.93,260) node [anchor=north west][inner sep=0.75pt]   [align=left] {$\displaystyle g_{a_{2}}^{-1} g_{a_{3}}$};
    \draw (468.93,260) node [anchor=north west][inner sep=0.75pt]   [align=left] {$\displaystyle g_{a_{m}}^{-1} g_{n}$};
    \draw (361.93,260) node [anchor=north west][inner sep=0.75pt]   [align=left] {$\displaystyle g_{a_{m-1}}^{-1} g_{a_{n}}$};
    \draw (204,108) node [anchor=north west][inner sep=0.75pt]   [align=left] {$\displaystyle s_{a_{i-1} +1} \cdots s_{a_{i}}$};
    \draw (29,276) node [anchor=north west][inner sep=0.75pt]   [align=left] {$\displaystyle 1$};
    \draw (524,279) node [anchor=north west][inner sep=0.75pt]   [align=left] {$\displaystyle \Phi _{g}( I)$};
    \draw (358,14) node [anchor=north west][inner sep=0.75pt]   [align=left] {$\displaystyle g_{a_{i+1}} o$};
    \draw (9,27) node [anchor=north west][inner sep=0.75pt]   [align=left] {$\displaystyle g_{a_{i} -L^{2}} o$};
    \draw (112,173) node [anchor=north west][inner sep=0.75pt]   [align=left] {$\displaystyle  h\cdot g_{a_{i} +L^{2}} o$};
    \draw (331,198) node [anchor=north west][inner sep=0.75pt]   [align=left] {$\displaystyle  h\cdot g_{a_{i+1}} o$};
    \draw (224,167) node [anchor=north west][inner sep=0.75pt]   [align=left] {$\displaystyle g_{a_{i}}^{-1} g_{a_{i} +L^{2}}$};

    \end{tikzpicture}
    
    \caption{The construction of $\Phi_g(I)$, where $h = (g_{a_i} f_i f) p_i (g_{a_i} f_i f)^{-1}\in H$}
    \label{fig:Con_Ext}
\end{figure}
\begin{remark}
    The existence of $f_i$ relies on the choice of $\mathcal{F}$ where we take the 5-tuple to be $$(g_{a_i}^{-1}g_{a_{i-1}}, g_{a_i}^{-1}g_{a_{i+1}}, g_{a_i}^{-1}g_{a_i-L^2}, g_{a_i}^{-1}g_{a_i+L^2}, f).$$

	The containment of the image of $\Phi_g$ in $Hg$ crucially uses property (6) of the confined subgroup $H$; cf. \cite[Lemma 4.6]{DY24:GrowthTight}. 
	
	If we identify using indicator function $\chi_I$ the subset $I$ of $\{1,\cdots,n\}$   as a sequence $(\epsilon_i)_{i=1}^n$ in $\{0,1\}^n$ by $\epsilon_i = \chi_I(i)$, then we can   write compactly $\Phi_g(I) = \prod_{i=1}^{n} s_i (g_i q_i g_i^{-1})^{\epsilon_i}$ where $g_{a_i} = f_i f$, $q_{a_i} = p_i$ for $1\le i\le m$ and $g_j = q_j = id$ if $j\neq a_i$ for any $1\le i\le m$.
\end{remark}
\begin{lem}\label{extCoset}
    For every $I\in \mathcal{P}(n, L)$, $\Phi_g(I)\in Hg$.
\end{lem}
\begin{proof}
    We prove by inducton on $m=|I|$. For the base case $m=0$, $I=\emptyset$ and $\Phi_g(I) = g\in Hg$.
    
    Now assume that $\Phi_g(I)\in Hg$ whenever $|I|=m-1$. Let $I\in\mathcal{P}(n, L)$ with $|I|=m$, and write $I = \{a_1,\cdots, a_m\}$ where $1\le a_1<a_2<\cdots<a_m\le n$. Let $J = I\setminus\{a_1\} = \{a_2, \cdots, a_m\}$. Then $|J|=m-1$ and we have $$\Phi_g(I) = g_{a_1} f_1 f p_1 f^{-1} f_i^{-1}g_{a_1}^{-1}\Phi_g(J) = (g_{a_1}f_1 f) p_1 (g_{a_1}f_1 f)^{-1} \Phi_g(J).$$ By induction hypothesis, $\Phi_g(J)\in Hg$. Furthermore, property (6) of \ref{defExtensionMap} implies $(g_{a_1}f_1 f) p_1 (g_{a_1}f_1 f)^{-1}\in H$. Thus, $\Phi_g(I)\in Hg$. This completes the induction step, and thus the lemma holds.
\end{proof}
\begin{lem}\label{extAdm}
	There exists $L_1>0$ such that $\gamma_g(I)$ is an $(D, B)$-admissible path for $I\in \mathcal{P}(n, L)$ with $L>L_1$. In particular, $\gamma_g(I)$ has the $r$-fellow travel property for some $r = r(D, B) > 0$.
\end{lem}
\begin{proof}
    Since the action of $G$ on $X$ is proper, let $L_1 = \max\{\|g\|_S: g\in G, d_X(o, go) \le B\}$.
    
	By the construction of $\gamma_g(I)$, it suffices to show that $\ax(f_i) \neq g_{a_{i}}^{-1} g_{a_{i+1}} \ax(f_{i+1})$.  
	
	If not, let $h = g_{a_{i}}^{-1} g_{a_{i+1}}$ and $Y = \ax(f_i) = h \ax(f_{i+1})$. Then $o, po\in Y$. So $\|h\|_X = d(o, h o) \le d^{\pi}_Y([o, po]) \le B$ by condition (\ref{nextproj}) in Definition \ref{defExtensionMap}. However, since $L>L_1$ and $\|h\|_S \ge L$, we have $\|h\|_X  > B$, which yields a contradiction.
\end{proof}

\subsection{Injectivity of Extension Map}
\label{subsec:injectivity}
In this subsection, we will prove that the extension map $\Phi_g$ introduced before is injective in $\mathcal{P}(n, L)$. Properties (1), (2), (3), (4) and (5) in Definition~\ref{defExtensionMap} will all be used to establish this, as we will see.

First, we need a preliminary result on the bounds of the projections onto two contracting subsets. This result implies that for any point, either one of the two projections lies near the common perpendicular of the two contracting sets, provided that the sets are far away from each other.

\begin{lem}\label{projToAxis}
	Let $Y, Z$ be $C$-contracting subsets with $B$-bounded projection. Let $\tau > 0$. Let $y\in Y$ and $z\in Z$ such that $d^\pi_Y(y, z)\le \tau$, $d^\pi_Z(y, z)\le \tau$ and $d(y, z) > 2\tau + 10C + B$. Then there exists $E = E(C, B)$ such that $$\min\{d^\pi_Y(x, y), d^{\pi}_Z(x, z)\} \le \tau + E$$ for any $x\in X$. 
\end{lem}
\begin{proof}
	Let $E = 2C + \hat C + B >C$. Suppose the opposite holds, that is $d^\pi_Y(x, y) > \tau + E$, $d^{\pi}_Z(x, z) > \tau + E$. Then $d^{\pi}_Z(x, y) \ge d^\pi_Z(x, z)- d^\pi_Z(y, z) > E > C$, which means that $[x, y]\cap N_C(Z)\neq\emptyset$. Noticing that $d(y, Z) \ge d(y, z) - d^\pi_Z(y, z) > E > C$, let $p$ be the first point on $[y, x]$ such that $d(p, Z)\le C$. Let $q \in \pi_Z(p)$ and $y' \in \pi_Z(y)$. Then by $C$-contracting property of $Z$, we have $d(q, y') = d^\pi_Z(y, p) \le C$. Since $Y, Z$ has $B$-bounded projection, $d^\pi_Y(q, z) \le B$. Considering $d(z, Y) \ge d(z, y) - d^{\pi}_Y(z, y) \ge \tau + 10C + B$, by triangular inequality we have $d(q, z) \le d(q, y') + d^\pi_Z(y, z) \le \tau + C$, $d(q,Y) \geq d(z,Y) - d(q,z) - d^\pi_Y(q,z)$, and hence $$d(q, Y) \ge (\tau + 10C + B) - (\tau + C) - B = 9C.$$ Thus, $$d(p, Y) \ge d(q, Y) - d(p, q) \ge 8C.$$ We claim that $[p, x]\cap N_C(Y) = \emptyset$.

    To prove the claim, assume, for the sake of contradiction, that $v$ is the first point on $[p,x]$ such that $d(v, Y) \le C$. Noticing that $[y, p]\cap N_C(Y)\neq \emptyset$, let $u$ be the last point on $[y, p]$ such that $d(u, Y) \le C$. Then by $C$-contracting property of $Y$, $d^\pi_Y(u, v)\le C$ and $d^\pi_Y(v, p)\le C$. So by triangular inequality we have $$d(u, v) \le d(u, Y) + d(v, Y) + d^\pi_Y(u, v) \le 3C.$$ 
    
    Now, since $y\in [u, v]$, we have $$d(u, v) \ge d(v, p) \ge d(p, Y) - d(v, Y) - d^\pi_Y(v, p) \ge 8C - C - C = 6C > 3C,$$ which is a contradiction.

    By the claim, we have $d^\pi_Y(p, x) \le C$. So, by triangular inequality we have $$d^\pi_Y(y, p)\ge d^\pi_Y(y,x) - d^\pi_Y(x,p) > \tau + E - C.$$ Then, by Lemma \ref{BigFive}, we get $d^\pi_Y(q, y) \ge d^\pi_Y(p, y) - d(p, q) - \hat{C} \ge \tau + E - 2C - \hat C$. Thus, by triangular inequality, it follows that $$d^\pi_Y(q, z)\ge d^\pi_Y(q, y)-d^\pi_Y(y,z) > E - 2C - \hat C.$$ However, recall that $q,z\in Z$ so $d^\pi_Y(q, z)\le B$ follows from $B$-bounded projection property. Since we took $E - 2C - \hat C = B$, we get a contradiction, and the lemma follows.
\end{proof}
Then we will prove the injectivity of $\Phi_g$.
\begin{lem}\label{extInj}
	Let $\pi: \Cay (G, S)\to X$ be a $(M, M)$-quasi-isometry and suppose $\pi(\gamma)$ is a $(K, C)$-quasi-geodesic for each word geodesic $\gamma\subset \Cay (G,S)$.
    Let $a=a(K,C)\ge 1$ be the constant in Lemma \ref{projLip} for $(K, C)$-quasi-geodesics and $C$-contracting subsets. There exists $L_2 = L_2(K, C, M, B)>0$ and $D = D(K, C, M, B) > 0$ so that the following holds.
	
	Let $g = s_1 \cdots s_n\in G$ be a word geodesic. For any $L>L_2$, if $\|f\|_X \ge 2a\log (L) + D$ for each $f\in \mathcal{F}$, then $\Phi_g: \mathcal{P}(n,L)\to Hg$ is injective.
\end{lem}
\begin{proof}
	Let $b, c$ be the constant in Lemma \ref{projLip} and let $E$ be in Lemma \ref{projToAxis}. Let $L_2>0$ be large enough such that $$LM^{-1} - M - 2(a\log (2ML^2) + b + cB) > 10 C + B$$ holds for each $L > L_2$, and $D$ be a large constant to be defined below. 
	
	Let $I_1, I_2 \in \mathcal{P}(n, L)$ with $I_1\ne I_2\in \mathcal{P}(n,L)$ and so that $h := \Phi_g(I_1) = \Phi_g(I_2)$. Let $I_1 = \{a_1, \cdots, a_{\alpha}\}$ and $I_2 = \{b_1, \cdots, b_{\beta}\}$. Let $f_i$ ($1\le i\le \alpha$) be the chosen contracting elements for $a_i\in I_1$ and $f_j'$ ($1\le j\le \beta$) be the chosen contracting elements for $b_j\in I_2$. We shall prove inductively that for each $1\le k\le \lfloor\frac{n}{L}\rfloor$, $$I_1 \cap \{L, 2L, \cdots, kL\} = I_2 \cap \{L, 2L, \cdots, kL\}$$ and simultaneously, $f_i = f_i'$ for each $1\le i\le \max\{\alpha,\beta\}$. This clearly    implies $I_1=I_2$ and $\alpha = \beta$, so the proof of injectivity follows.
	
	We start by proving  the base step   $k=1$. The following claim is the main observation in the base step and the subsequent induction:
	\begin{claim}
		$L\in I_1 \cap I_2$ or $L\not\in I_1\cup I_2$.
	\end{claim}
	\begin{proof}[Proof of the claim]
		If the claim is false, without loss of generality, we assume $L\in I_1$ and $L\not\in I_2$. Then $a_1 = L$ and $2L\le b_1 \le L^2$ since $L^2\in I_2$. To derive a contradiction, the idea is to first give a lower bound on  $d^{\pi}_Y(o, h o)$ using  the admissible path   $\gamma_g(I_1)$, and then a upper bound using     $\gamma_g(I_2)$ where $Y$ is a $C$-contracting subset of $X$.
		Recall that $g_i = s_1 s_2\cdots s_i$ for each $1\le i\le n$ and $g_0 = id$. Denote $Y = g_L \ax(f_1)$ be the $C$-contracting subset associated with $[g_L o, g_L f_1 o]$ in the admissible path   $\gamma_g(I_1)$.   By the $r$-fellow travel property, the two endpoints of  $[g_L o, g_L f_1 o]$ are contained in the $r$-neighborhood of   $[o, ho]$, so there exist  $u, v\in [o, ho]$ such that $d(u, g_L o)\le r$, $d(v, g_L f_1 o)\le r$. If  $u' = \pi_Y(u)$ and $v' = \pi_Y(v)$, the triangular inequality shows $$d(u', g_L o)\le d(u', u) + d(u, g_L o) = d(u, Y) + d(u, g_L o)\le 2 d(u, g_L o) \le 2r.$$ Similarly, $d(v', g_L f_1 o)\le 2r$. Thus, we obtain   $$d^{\pi}_Y(u, v) \ge d(u', v') \ge d(g_L o, g_L f_1 o) - d(u', g_L o) - d(v', g_L f_1 o) \ge \|f_1\|_X - 4r$$ which yields   by Lemma \ref{BigFive}  \begin{equation}\label{lowerProj}
			d^{\pi}_Y(o, h o) \ge \diam{\pi_Y([o, ho])} - \hat C \ge d^{\pi}_Y(u, v) - \hat C \ge \|f_1\|_X - 4r - \hat C.
		\end{equation}
		
		For a contradiction, we  use the other   path   $\gamma_g(I_2)$ to derive an upper bound for $d^{\pi}_Y(o, h o)$. Let $t:=b_1\le L^2 $ be the minimal integer in $I_2$. Our idea is to give an lower bound on $d(g_t o, Y)$ then use Lemma \ref{projToAxis} to give an upper bound on the projection. By construction of $f_1$ as in Definition \ref{defExtensionMap}(\ref{nextLproj}), we have $$d^{\pi}_Y(g_L o, g_{L^2+L} o) \le B.$$
		 Since $\pi$ is an $(M, M)$-quasi-isometry, we have $M^{-1}(j-i) - M\le d(g_i o, g_j o)\le M(j-i) + M$. According to   Lemma \ref{projLip}, since $\pi([g_L, g_{L+L^2}]_S)$ is a $(K, C)$-quasi-geodesic containing $g_{b_1}o=g_to$ and $d(g_L o, g_{L+L^2} o)\le M L^2 + M$, we have $$d^{\pi}_Y (g_L o, g_{t} o) \le a\log(ML^2 + M) + b + cB \le a\log(2ML^2) + b + cB.$$ 
        
        On the other hand, denote $Z = g_t\ax(f_1')$ be the $C$-contracting subset associated with $[g_t o, g_t f_1' o]$ in the admissible path $\gamma_g(I_2)$. Then by the same reason as above, noticing that by Definition \ref{defExtensionMap}(\ref{prevLproj}) $d^\pi_Z (o, g_t o)\le B$, we have $$d^{\pi}_Z (g_L o, g_{t} o) \le a\log(2ML^2) + b + cB.$$
        
		From here, we are going to apply Lemma \ref{projToAxis}. Let $\tau = a\log(2ML^2) + b + cB$, $y = g_L o$, $z = g_t o$. Then by the equalities above we have $d^\pi_Y(y, z)\le \tau$ and $d^\pi_Z(y, z)\le \tau$. By choice of $L_2$ and $L > L_2$ we have $d(y, z) \ge LM^{-1} - M > 2\tau + 10C + B$. Clearly $Y\neq Z$ for $d^\pi(y, z) < d(y, z)$ and thus $Y, Z$ have $B$-bounded projection. Since $\gamma_g(I_2)$ has $r$-fellow travel property, we have $d^\pi_Z(z, h o) \ge \|f_1'\| - r > \tau + E$. Suppose we have $D > \tau - 2a\log(L) + E + r = a\log(2M) + b + cB + E + r$. Then by Lemma \ref{projToAxis}, $d^{\pi}_Y(h o, g_L o) \le \tau + E$.
        
	    Also, property (\ref{prevproj}) in Definition \ref{defExtensionMap} implies $d^{\pi}_Y (o, g_L o)\le B$. Hence, the triangular inequality for $d^{\pi}_Y$ implies 
		\begin{equation}\label{upperProj}
			d^{\pi}_Y(o, h o)\le d^{\pi}_Y(o, g_L o) + d^{\pi}_Y(h o, g_L o)\le \tau + E + B.
		\end{equation}
		
		Finally, combining the inequalities (\ref{lowerProj}) and (\ref{upperProj}) we have $$\|f_1\|_X - 4r - \hat C \le d^{\pi}_Y(o, h o) \le \tau + E + B.$$ If we set $D = 4r + \hat C + a \log (2 M) + b + cB + E + B > \tau - 2a\log(L) + E + r$, this contradicts our assumption that $\|f_1\|_X > 2a\log L + D$. So the claim holds.
	\end{proof}
	
	If $1\in I_1\cap I_2$, it remains to prove $f_1=f_1'$ for the base case.  Indeed, if  $f_1\neq f_1'$, then by the choice of $\mathcal{F}$, $g_L \ax(f_1)$ and $g_L \ax(f_1')$ has bounded projection property. For a path $\gamma$, denote the subpath of $\gamma$ with endpoints $x, y\in \gamma$ by $\gamma|_{[x, y]}$. Let $\gamma' = \gamma_g(P_2)|_{[h o, g_L o]} \cup \gamma_g(P_1)|_{[g_L o, h o]}$. Then $\gamma'$ is a $(D, B)$-admissible path. Thus it is a $\lambda$-quasi geodesic loop at $h o$ so $\mathrm{length}(\gamma') \le \lambda^2$. This is in contradiction with $\mathrm{length}(\gamma') \ge \mathrm{length}(f_1) > F > \lambda^2$. So $f_1 = f_1'$. This finishes the proof of the base case $k = 1$.
	
	Let us consider $k \ge 2$. By   the induction hypothesis,  $I:= I_1 \cap \{L, \cdots, (k-1)L\} = I_2 \cap \{L, \cdots, (k-1)L\}$ and $f_i = f_i'$ for $a_i\le k-1$. Let $x = g_{(k-1)L}$, $y = g_{(k-1)L}^{-1} g$, $h_1 = \Phi_{x}(P)$ and $J_i = (I_i - (k-1)L)\cap \mathbb{N} = \{j: j\ge 1, j+(k-1)L\in I_i\}$ for $i = 1, 2$. Then we have $\Phi_g(I_i) = \Phi_x(I) \Phi_y(J_i)$ for $i = 1, 2$. Thus, $\Phi_y(J_1) = \Phi_y(J_2)$. By the case of $k=L$ we have $J_1 \cap \{L\} = J_2 \cap \{L\}$ and the corresponding contracting element is the same. By the definition of $J_i$ we can conclude that $I_1\cap \{L, \cdots, kL\} = I_2\cap \{L, \cdots, kL\}$ and $f_i = f'_i$ if $a_i\le k$.
	
	By   induction, we   conclude that $P_1 = P_2$ and $f_i = f_i'$ for each $i$. Thus the lemma holds.
\end{proof}
\begin{remark}
    In Lemma \ref{extInj}, $\pi, K, C, a$ are independent of the confined subgroup $H$. On the other hand, $B, L_2, D$ depend on $H$.
\end{remark}

\subsection{A Gap Criterion}
\label{subsec:gap}
In this section we prove the ``gap criterion'', giving sufficient conditions on a generating set $S$ and subset $A \subset G$ to satisfy $\omega_A(d_{S}) < \omega_G(d_{S})$.

The proof uses the critical exponent series of a subset $A\subset G$, which we define now.
Given $g \in G$, let $\|g\|_S$ denote the \emph{word norm} of $g$ with respect to the generating set $S$. For $A \subset G$ and a real parameter $s \geq 0$, define the \emph{Poincar\'e series} of $A$ as
    \begin{equation}
        P_A(s) = \sum_{g\in A} e^{-s\|g\|_S},
        \label{eq:critical_exp_series}
    \end{equation}
    which converges for $s > \omega_A(d_S)$ and diverges for $s < \omega_A(d_S)$.

Loosely speaking, for a subset $A \subset G$, the gap criterion provides conditions under which one can construct, for each element $g\in A$, exponentially many distinct elements in $G$ of word-length roughly $\| g\|_S$.
Moreover, for any $g' \in G$ distinct from $g$, $\mathrm{Im}(\Phi_g) \cap \mathrm{Im}(\Phi_g')$ are disjoint. Since the critical exponent captures exponential growth rates, this yields $\omega_A(S) < \omega_G(S)$.

Fix  a finite generating set  $S$ of $G$ and denote $\omega_S = \omega_G(d_{S})$.  
A similar result appears in \cite[Lemma 2.26]{DY24:GrowthTight}.

\begin{lem}\label{gapCritLem}
	Let $\theta, M>0$ and $0<a \le \frac{1}{8\omega_S}$. Denote  $L_3 = 4 e^{8\omega_S M + 2} + 2$.  Let $L > L_3$ and  $A\subset G$ be a subset  so that the following holds:
	\begin{enumerate}
		\item \label{Conditiondelta} For any $g\in A$, there exists an injective map $\Phi_g: \mathcal{P}(n, L) \to G$ where $n = \|g\|_S$;
		\item \label{ConditionL} For any $g\in A$ and $I\in \mathcal{P}(n, L)$, we have $\|\Phi_g(I)\|_S \le \|g\|_S + |I| (a  \log(L) + M)$;
		\item \label{ConditionD} For any $g_1 \neq  g_2\in A$, $\mathrm{Im}(\Phi_{g_1}) \cap \mathrm{Im}(\Phi_{g_2}) = \emptyset$.
	\end{enumerate}
	Then $\omega_A(d_{S}) \le \omega_S - \frac{1}{L^2}$.
\end{lem}

We note that, in applying this result, the hypothesis that $g \in A$ in condition (1) of Lemma~\ref{gapCritLem} will be superfluous, since in Sections~\ref{subsec:extensionmap} and~\ref{subsec:injectivity} we have constructed a map $\Phi_g$ satisfying condition (1) of Lemma~\ref{gapCritLem} \emph{for every} $g \in G$. However, the proof of Lemma~\ref{gapCritLem} will only use that the map $\Phi_g$ is defined for $g \in A$.  
In applications, the choice of $g \in A$ will play a role in conditions (2) and (3).
\begin{proof}
Recall that the critical exponent series of $A$ it is defined as $$P_A(s) = \sum_{g\in A} e^{-s\|g\|_S},$$
	it converges at $s>\omega_A(d_{S})$ and diverges at $s<\omega_A(d_{S})$.
	For any $s>\omega_S$, we have, 
	\begin{equation}\label{ineqPS}
		\begin{aligned}
			P_G(s) &= \sum_{g\in G}  e^{-s\|g\|_S} & \\
			&\ge \sum_{g\in A} \sum_{h\in \mathrm{Im}(\Phi_g)} e^{-s \|h\|_S}&\mbox{by }(\ref{ConditionD})\\
			&\ge \sum_{g\in A} \sum_{I\in \mathcal{P}(n, L)} \exp\left(-s (\|g\|_S + |I| (a \log(L) + M))\right) & \mbox{by }(\ref{ConditionL})\\
			&= \sum_{g\in A} e^{-s \|g\|_S} \sum_{I\in \mathcal{P}(n, L)} \exp(-s |I| (a \log(L) + M)))&\\
			&= \sum_{g\in A} e^{-s \|g\|_S} \sum_{I\in \mathcal{P}(n, L)} (L^{-sa} e^{-sM})^{|I|}&
		\end{aligned}
	\end{equation}
	Denote $F(g) := \sum_{I\in \mathcal{P}(n, L)} (L^{-sa} e^{-sM})^{|I|}$. Let $\epsilon_i = \chi_I(i)$ for $1\le i\le n$, $m = \lfloor\frac{n}{L}\rfloor$ and by Binomial Theorem we have 
	\[\begin{aligned}
		F(g) &= \sum_{(\epsilon_1, \cdots, \epsilon_n)\in \mathcal{P}(n, L)} (L^{-sa} e^{-sM})^{\epsilon_1 + \cdots + \epsilon_n}\\
		&= (L^{-sa} e^{-sM})^{\lfloor\frac{n}{L^2}\rfloor} (1+L^{-sa} e^{-sM})^{[\frac{n}{L}] - [\frac{n}{L^2}]}\\
		&\ge (L^{-sa} e^{-sM})^{\frac{m}{L}} (1+L^{-sa} e^{-sM})^{m - \frac{m}{L}}\\
		&= \left((L^{-sa} e^{-sM}) (1+L^{-sa} e^{-sM})^{L-1}\right)^{\frac{m}{L}}.
	\end{aligned}\]
	where the second equality holds since $\epsilon_i$ must be 1 if $i = kL^2$ and $\epsilon_i$ can be either 0 or 1 otherwise.\\
	Fix any $s$ such that $\omega_S < s < 2\omega_S$. Then $sa < \frac{1}{4}$ and we have the following inequality:
	\[\begin{aligned}
		F(g)^{\frac{L}{m}} &\ge (L^{-sa} e^{-sM}) (1+L^{-sa} e^{-sM})^{L-1} &\\
		&\ge L^{-\frac{1}{4}} e^{-sM} (1+L^{-\frac{1}{4}} e^{-sM})^{L-1} &  \because (1+x)^n  > nx \\
		&\ge L^{-\frac{1}{4}} e^{-sM} L^{-\frac{1}{4}} e^{-sM}(L-1) & \because  L>2\\
		&\ge L^{-\frac{1}{4}} e^{-sM} L^{-\frac{1}{4}} e^{-sM}\frac{L}{2} = \frac{1}{2} e^{-2sM} L^{\frac{1}{2}} \ge e &
	\end{aligned}\]
	where the last inequality holds since $L^{\frac{1}{2}} > 2 e^{4\omega_S M + 1} \ge 2 e^{2s M + 1}$. Combined with condition (1), $n = \|g\|_S$, we obtain \[F(g) \ge e^{\frac{m}{L}} \ge e^{\frac{n}{L^2} - \frac{1}{L}}.\]
	Plugging the above inequality into (\ref{ineqPS}), we have
	\[\begin{aligned}
		P_G(s) &\ge \sum_{g\in A} e^{-s \|g\|_S} e^{\frac{ \|g\|_S}{L^2} - \frac{1}{L}} = e^{-\frac{1}{L}} \sum_{g\in A} e^{-(s - \frac{1}{L^2}) \|g\|_S} = e^{-\frac{1}{L}} P_A(s - \frac{1}{L^2}).
	\end{aligned}\]
	Notice that $s > \omega_S$ implies $P_G(s)$ converges and $\omega_S < s < 2\omega_S$ can be chosen arbitrarily. So $P_A(s - \frac{1}{L^2})$ also converges for any $2\omega_S > s > \omega_S$. Thus, $\omega_A(d_{S}) \le \omega_S - \frac{1}{L^2}$.
\end{proof}
\subsection{Approximating by Generators}
Our goal in this section is to construct the generating sets approximating $X$, and to show they satisfy the gap criterion, thus proving Theorem~\ref{GrowthTight4SomeS}.

Let $\omega_X = \omega_G(d_X)$ be the growth rate of the action. The following result is just a restatement of Lemma~\ref{lem:density}.

\begin{lem}\label{approByBall}
	Suppose $G$ acts properly and co-boundedly   on $X$. Let $n$ be large enough. Then $S = B(o, n)$ generates $G$ and there exists $N>0$ so that $\frac{n}{2} > N > 0$ and $(n-N) (d_S(x, y)-1) \le d_X(x, y) \le n d_S(x, y)$. Furthermore, $(n - N) \omega_X \le \omega_G(d_{S}) \le n \omega_X$.
\end{lem}

We now proceed with the proof of Theorem~\ref{mainthm:cogrowth}. The main challenge lies in constructing a generating set that is independent of the confined subgroup $H$. In our approach, we define the generating set as a sufficiently large ball, ensuring that the quasi-geodesic constants $(K, C)$ remain bounded, and in particular $K$ is close to 1. This, in turn, guarantees that the constant $a$ in Lemma~\ref{projLip} is sufficiently small, depending only on $\omega_X$. 

For a given confined subgroup $H$, we utilize the extension map constructed in Definition~\ref{defExtensionMap}. To ensure injectivity in Lemma~\ref{extInj}, we require the $C$-contracting axes to be sufficiently long—specifically, of length greater than $2a \log(L) + D$, where $a$ is independent of $H$ and $D$ depends on $H$. Finally, we apply Lemma~\ref{gapCritLem} to establish a growth gap.

For convenience, we restate the result here.

\begin{thm}\label{thm:main_growth}
Let $G$ act properly coboundedly on a geodesic space $(X,d)$ with a strongly contracting element, and let $x_0 \in X$. Then, there exists $n_0 > 0$ such that for all $n \geq n_0$, if $S_n$ consists of elements whose $x_0$-orbit lies in $B(x_0,n) \subset X$, then $S_n$ generates $G$ and satisfies 
\[
\omega_{G/H}(d_{\bar{S}}) < \omega_G(d_S)
\]
for \emph{every} infinite confined subgroup $H \subset G$ with a non-degenerate confining subset.
\end{thm}
\begin{proof}
    Let $n_0 \in \mathbb{N}$ be chosen as in Construction~\ref{cons:generatingsets}, so that any of the sets $S_n$ for $n \geq n_0$ generates $G$.
    To lighten notation, we will  write $S \coloneqq S_n$.
	By Lemma \ref{approByBall}, $\pi([x, y]_S)$ is an $(\frac{n}{n-N}, n)$-quasi-geodesic for any $x,y\in G$.  Let $a, b, c$ be the corresponding constants in Lemma \ref{projLip}. Since $\lim_{n\to \infty} \frac{n}{n-N} = 1$, we may choose    some large $n > 2N$ so that   $a < \frac{1}{64\omega_X}$. Moreover, for any $C$-contracting subset $Y$ and $(\frac{n}{n-N}, n)$-quasi-geodesic $\gamma$ with an endpoint in $Y$ we have $\diam{\pi_Y(\gamma)} \le a \log(\mathrm{length}(\gamma)) + b + c \cdot d^{\pi}_Y(\gamma_-, \gamma_+)$.
	
	Let $H$ be a confined subgroup with a non-degenerate confined subset $P$. Let $A$ be a set of shortest representatives of $G/H$: for any $g\in A$, $\|g\|_S=d_S(1,Hg)$. Thus, $\omega_{G/H}(\bar{S}) = \omega_A(d_{S})$. 
	
	Let $L_1$, $L_2$, $D$ be the constants in Lemmas \ref{extAdm} and \ref{extInj}.
	
	Let $M = 2(\frac{4D}{n} + 1) + \max\{\|p\|_S:p\in P\} + 2\|f\|_S$, which depends on $P$ (and $H$). 
	Let $L_3 =L_3(M)$ be given by Lemma \ref{gapCritLem}. Choose $L > \max\{L_1, L_2, L_3\}$. To conclude  $\omega_{G/H}(d_{\bar{S}}) = \omega_A(d_{S}) < \omega_G(d_{S})$, it suffices to verify the conditions (\ref{Conditiondelta}), (\ref{ConditionL}), (\ref{ConditionD}) in Lemma \ref{gapCritLem}.
	
	By taking a sufficiently large power, we may assume that $2a \log(L) + D \le \|f_i\| \le 2a \log(L) + 2D$ for each $f_i\in \mathcal{F}$.
	Thus, for any $g\in A$  the extension map $\Phi_g: \mathcal{P}(n, L) \to Hg$  defined in Definition \ref{defExtensionMap} is injective by Lemma \ref{extInj}. This furnishes the injective map in (\ref{Conditiondelta}). Since $\mathrm{Im}(\Phi_g) \subset Hg$ and $A$ consists of representatives of $G/H$, we have $\mathrm{Im}(\Phi_{g_1})\cap \mathrm{Im}(\Phi_{g_2}) = \emptyset$ for any $g_1\ne g_2\in A$, so  (\ref{ConditionD}) holds. It remains to prove (\ref{ConditionL}) in what follows. 
	
	By Lemma \ref{approByBall}, we have $$\forall f\in \mathcal F,\quad \|f\|_S \le \frac{1}{n-N} \|f\|_X + 1 \le \frac{2}{n}(2a \log(L) + 2D) + 1.$$
	which yields by the choice of $M$ $$\|f_i f p_i f^{-1} f_i^{-1}\|_S \le 2\|f_i\|_S + 2\|f\|_S + \|p\|_S \le \frac{8a}{n} \log(L) + M.$$ 
	
	Notice that $$\frac{8a}{n} < \frac{1}{8\omega_X n} \le \frac{1}{8\omega_G(d_{S})}.$$
	
	Thus, the condition (\ref{ConditionL}) is fulfilled: $$\|\Phi_g(I)\|_S\le \|g\|_S + |I| \left(\frac{1}{8\omega_G(d_{S})} \log(L) + M\right).$$

	In conclusion, we have verified all the conditions  in Lemma \ref{gapCritLem}, so it follows that $\omega_{G/H}(d_{\bar{S}}) = \omega_A(d_{S}) < \omega_G(d_{S})$. The proof is complete.
\end{proof}

\begin{remark}
Our arguments do not specifically rely on the fact that $d_S$ is a word metric; rather, they apply to any roughly geodesic metric. In particular, if one can construct a sequence of roughly geodesic metrics $d_n$ that approximates a metric with strongly contracting elements, then growth tightness holds for all $d_n$ with $n \geq n_0$. 

For example, in forthcoming work by the second author~\cite{CMGR:GreenDensity}, Green metrics arising from admissible random walks on $G$ are shown to be dense when $G$ is a non-elementary Gromov hyperbolic group. One might hope to extend this result beyond Gromov hyperbolic groups, at least for injective metrics. Specifically, can every injective metric in the mapping class group be approximated by a sequence of Green metrics?
\end{remark}

\appendix

\section{Many injective metrics on mapping class groups}
\label{subsec:proof_infinitely}

In this Appendix we prove Theorem~\ref{thm:infinitely_many_injective}.
 The argument for constructing even a single injective metric on the mapping class group is quite intricate. However, this task has already been accomplished in \cite{haettelhodapetyt:coarse} and \cite{PetytZalloum24} through different methods. To keep this section concise, we provide only a sketch of the proof of Theorem \ref{thm:infinitely_many_injective}, directing the interested reader to \cite{PetytZalloum24} for further details.

 For the rest of this section, we will fix $G$ to be the mapping class group of a finite type surface $\Sigma$. Unlike in the other results of this paper, this restriction will be important since we will be using that $\MCG(\Sigma)$ is acylindrically hyperbolic.
 
 We will often think of $G$ as the Cayley graph resulting from a fixed once and for all finite generating set. When we say that something depends only on $G$, we mean that it also depends on this choice of a finite generating set. For each subsurface $Y \subset \Sigma,$ we have a (coarse) map $\pi_Y:G \rightarrow \mathcal{C}(Y)$ which is $(E,E)$-coarsely Lipschitz. Here, $\mathcal{C}(Y)$ denotes the \emph{curve graph} of the subsurface $Y.$ The constant $E$ depends only on $G$ and not on the surface $Y.$

The core idea behind constructing the injective metrics $\{\sigma_i\}_{i \in \mathbb{N}}$ in Theorem \ref{thm:infinitely_many_injective} is straightforward: replace the metric space $(\mathcal{C}(\Sigma),d)$ with $(\mathcal{C}(\Sigma),\frac{d}{i})$ and follow the construction of the injective metric in \cite{PetytZalloum24} for the new HHS structure on $G$. A key observation is that the resulting injective metrics on $G$ are \emph{not} scalar multiples of each other, as only $\mathcal{C}(\Sigma)$ is rescaled, while $\mathcal{C}(Y)$ remains unchanged for all proper subsurfaces $Y \subset \Sigma$. For an intuitive example, consider $\mathbb{R}^2$ with its standard $\ell^\infty$-metric versus the $\ell^\infty$-metric obtained by modifying the metric on the second $\mathbb{R}$ factor to $\frac{d_2}{n}$.

\subsection{Coarse Medians and Convex Hulls}
We begin by recalling some convexity notions in the mapping class group and related facts. For more details, see \cite{bowditch:convex}. 

A \emph{coarse median space} $(X,d, \mu)$ is a metric space $(X,d)$ equipped with a ternary operator $\mu$ (called the \emph{coarse median}), such that every finite subset can be approximated by a finite $\CAT(0)$ cube complex, with distortion controlled by the metric. The mapping class group $G$ admits a coarse median: for any $x,y,z \in G$, there is a (coarsely-defined) point $\mu(x,y,z)$. 

A subset $A \subset G$ is said to be \emph{$K$-convex} if for any $x,y \in A$ and $z \in G$, we have $d(\mu(x,y,z), A) \leq K.$ Given a subset $A \subset G$, one can associate a $K$-convex set, denoted $\hull(A)$, that contains $A$. The constant $K$ depends only on $G$ and is independent of $A$. Given a $K$-convex set $A$, there exists a natural (coarse) map $G_A: G \to A$, which is $(K',K')$-coarsely Lipschitz. The coarseness of this map, as well as the constant $K'$, depends only on $K$ and $G$, ensuring uniformity when $K$ is fixed.

\begin{remark}\label{rmk:constant_track} Given $A \subset G,$ there are various equivalent ways to construct $\hull(A)$. One way (described in \cite{Russell2018ConvexityIH}) is to project $A$ to each $\mathcal{C}(Y)$, connect pairs of points in $\pi_Y(A)$ by geodesics, and take all the points in $G$ whose projection lies  $E$-close to $\pi_Y(A).$ In particular, if $A \subset G$ is such that $\pi_Y(A)$ is $E$-convex for in some $\mathcal{C}(Y)$, then $\pi_Y(\hull(A)) \subset N(\pi_Y(A), 2E)$.
    
\end{remark}

\subsection{Weighted Curtains}\label{sec:fixed_p.A} For the rest of the section, we will fix a pseudo-Anosov $g \in G.$ Since $g$ is a pseudo-Anosov, it must act loxodromically on $\mathcal{C}(\Sigma)$, and there is a constant $D(g)>0$ such that for any proper subsurface $Y \subsetneq \Sigma,$ we have $\text{diam}_{\mathcal{C}(Y)}(\pi_Y(\langle g \rangle)) \leq D(g)$. Our ultimate goal is to construct a sequence of injective metrics $\rho_n \in \DDG$ where $\ell_{\rho_n}(g) \rightarrow 0$.

As mentioned in Subsection \ref{subsec:construction_injective}, the primary idea behind constructing an injective metric on $G$ is to define the distance between two points $x,y \in G$ using their ``amended" distance formula sum. Here, ``amended" means replacing the partial sum from disjoint surfaces with the maximal distance attained in a subsurface from that collection. To formalize this construction, we employ the notion of a \emph{curtain} (see \cite{PSZCAT}, \cite{zalloum:injectivity}, and \cite{PetytZalloum24}).

\begin{defn}[Weighted Curtains] \label{def:weights}
Fix $n \in \mathbb{N}$. 
\begin{enumerate}
    \item Let $Y \subsetneq \Sigma$. A \emph{$Y$-curtain} is defined as:
    \[
    h_{Y,b,I}:=\hull(\pi_Y^{-1}\pi_b^{-1}(I))
    \]
    where $b \subset \mathcal{C}(Y)$ is a geodesic, $I \subset b$ is a subgeodesic of length $D(g)+10E$, not containing any endpoints of $b$, and $\pi_b:U \to b$ is the (coarse) nearest-point projection map.
    
    \item A \emph{$\Sigma(n)$-curtain} is defined as:
    \[
    h_{\Sigma,b,I}:=\hull(\pi_\Sigma^{-1}\pi_b^{-1}(I))
    \]
    where $b \subset \mathcal{C}(\Sigma)$ is a geodesic, $I \subset b$ is a subgeodesic of length $10E n$, not containing any endpoints of $b$, and $\pi_b: \mathcal{C}(\Sigma) \to b$ is the (coarse) nearest-point projection map.
\end{enumerate}
\end{defn}

Let $\mathcal{H}^n$ be the collection of all $Y$-curtains and $\Sigma(n)$-curtains (with $n$ fixed throughout this discussion). A \emph{curtain} is defined as an element of $\mathcal{H}^n$. Despite the asymmetry of weights in the definition of curtains, the invariance of $\mathcal{C}(\Sigma)$ under $G$ ensures that our collection of curtains remains invariant under $G$. Further, it is easy to check that in each $\mathcal{C}(Y)$ (where $Y$ is allowed to be $\Sigma)$; the set $\pi_b^{-1}(I)$ is $E$-convex. In particular, via Remark \ref{rmk:constant_track}, the projection of $\hull(\pi_Y^{-1}\pi_b^{-1}(I))$ to $\mathcal{C}(Y)$ is $2E$-close to $\pi_b^{-1}(I)$.

\begin{defn}[Half-Spaces and Chains] Let $h \in \mathcal{H}^n$ and let $I,b,Y$ be as in Definition \ref{def:weights}. Here, $Y$ is allowed to be $\Sigma.$ The set $b-I$ consists of two connected components $I_1, I_2$. The \emph{half spaces} corresponding to $h$ are defined to be the sets $h^+=\hull(\pi_Y^{-1}\pi_b^{-1}(I_1))$ and $h^-= \hull(\pi_Y^{-1}\pi_b^{-1}(I_2))$. Now, if $d$ denotes the word metric on $G$, the fact that $\pi_Y$ is $(E,E)$-coarsely Lipschitz combined with Remark \ref{rmk:constant_track} give us that $d(h^+,h^-) \geq E$. Further, by definition, $h^+,h^-$ are $K$-convex (for a uniform $K$ that's independent of $n)$ and satisfy $G-h \subset h^+ \cup h^-$. A curtain $h$ is said to \emph{separate} subsets $A,B \subset X$ if $A \subset h^+$ and $B \subset h^-$. A collection of curtains $\{h_1, \dots, h_k\}$ is called a \emph{chain} if each $h_i$ separates $h_{i-1}$ from $h_{i+1}$. 
\end{defn}

\subsection{Building the Injective Metrics} 
\label{subsec:coarsely_injective}
A metric space is said to be \emph{coarsely injective} if there exists a constant $\delta > 0$ such that for any family $\{ B(x_i, r_i) : i \in I \}$ of balls satisfying 
\[
d(x_i, x_j) \leq r_i + r_j \quad \text{for all } i, j \in I,
\]
the $\delta$-neighborhood of these balls has a non-empty total intersection. This notion is a coarse analogue of injective metric spaces. A metric space is \emph{injective} if for any family $\{ B(x_i, r_i) : i \in I \}$ of balls satisfying 
\[
d(x_i, x_j) \leq r_i + r_j \quad \text{for all } i, j \in I,
\]
the total intersection of these balls is non-empty. This is equivalent to the definition of injective space we have given in Section~\ref{subsec:construction_injective}, by~\cite[Proposition~2.3]{Lang13:Injective}.

The \emph{injective hull} of a metric space $X$ is an injective metric space $E(X)$, together with an isometric embedding
\[
e \colon X \to E(X),
\]
such that no proper injective subspace of $E(X)$ contains $e(X)$ (see \cite[Section~3]{L13:Injective}).

A subset $Y$ of a metric space $X$ is \emph{coarsely dense} if there exists $r > 0$ such that every $x \in X$ is within distance $r$ of some $y \in Y$. The notions of injectivity and coarse injectivity are related through the injective hull: a metric space is coarsely injective if and only if it is coarsely dense in its injective hull \cite[Proposition~1]{haettelhodapetyt:coarse}.

Recall that $G$ is fixed to be a mapping class group of a finite-type surface $\Sigma.$ We now construct a collection of $G$-invariant coarsely injective metrics on $G$ each of which is quasi-isometric to any word metric on $G$ (though with varying quasi-isometry constants). These injective spaces will be used in the proof of Theorem \ref{thm:infinitely_many_injective} to establish that the space $\DDGS$ has infinite diameter.

\begin{defn}[$n$-Injective Metrics] For $x,y \in G$ with $x=y,$ define $\rho_n(x,y)=0$. If $x \neq y,$ define 
\[
\rho_n(x,y) := 1 + \max\{|c| : c \text{ is a chain separating } x \text{ and } y\}.
\]
\end{defn}

It is straightforward to verify that $\rho_n$ defines a metric. Moreover, Lemma 7.2 (or Theorem 0.1 in the Appendix) of \cite{zalloum:injectivity} establishes that this metric is quasi-isometric to the original word metric $d$ on $G$. However, it is important to note that the quasi-isometry constants depend heavily on the parameter $n$.

The proof that the above metric is injective is essentially contained in \cite{PetytZalloum24}, where the authors show that any metric defined via a sufficiently nice collection of walls is coarsely injective. Below, we provide a description of the argument. Before we do that, we shall fix some constants.

As discussed earlier, the coarse maps $\pi_Y:G \rightarrow \mathcal{C}(Y)$ are $(E,E)$-coarsely Lipschitz for a uniform $E.$ For any set $A \subset G$, the set $\hull(A)$ is $K$-convex for a uniform $K.$ From now and on, we will assume that $E$ is taken large enough so that every $\hull(A)$ is $E$-convex (in other words, we have just increased the coarse constant for the maps $\pi_Y).$ For a subset $A \subset G$, we use $N(A,K)$ to denote points $x \in G$ such that $d(x,a) \leq K$ for some $a \in A.$

\begin{lem} 
There exists a constant $\alpha$, depending only on $G$ (and not on $n$), such that every pair of points in $G$ is connected by an $\alpha$-rough geodesic for the metric $\sigma_n.$
\label{lem:rhon_geodesic}
\end{lem}

\begin{proof} 
Let $x, y \in G$ and fix an integer $r > 0$. Define $\mathcal{H}^r(x)$ to be the collection of curtains $h \in \mathcal{H}^n$ such that there exists a chain of length at least $r+1$ separating $x$ and $h$. Define 
\[
J_n(x,r) := \bigcap_{h \in \mathcal{H}^r(x), \, x \in h^+} N(h^+,E).
\]
Since each $h^+$ is  $E$-convex, the intersection $J_n(x,r)$ is also $E$-convex. We construct the rough geodesic as follows: consider the sequence $J_n(x,1), J_n(x,2), \dots, J_n(x,\rho_n(x,y)).$ Since each $J_n(x,r)$ is $E$-convex, there exists a uniformly coarsely well-defined gate map $G_r: G \to J_n(x,r)$. The sequence $G_1(y), G_2(y), \dots, G_{\rho_n(x,y)}$ provides a uniformly controlled rough geodesic. Verifying that such a sequence is a rough geodesic is a technical argument which is very similar to the one given in \cite[Proposition~4.6]{PetytZalloum24} (see also Figure 1 in \cite[Figure~1]{PetytZalloum24}). For compactness of the exposition, we omit the argument. 
\end{proof}

The following lemma establishes that balls with respect to the metric $\rho_n$ are $A$-convex for a uniform constant $A$ depending only on $G$ (and not on $n$). The proof follows the same argument as Lemma 4.2 in \cite{PetytZalloum24}, but we provide it here for completeness.

\begin{lem} 
There exists a constant $C$ depending only on $E$ such that for any ball $B_n(x,r) \subset (X,\rho_n)$, we have
\[
d_{\text{Haus}}(B_n(x,r), J_n(x,r)) < C.
\]
\end{lem}

\begin{proof} 
Suppose $y \in B_n(x,r)$ but $y \notin J_n(x,r)$, that is,  $y \notin N(h^+,E)$ for some $h \in \mathcal{H}^r(x)$. By definition of $h$, there exists a chain $c$ of length at least $r+1$ separating $x$ from $h$, and $c$ separates $x$ from $h \cup h^-$. In particular, $c$ also separates $y$ from $x$, contradicting the assumption that $y \in B_n(x,r)$. This establishes that $B_n(x,r) \subset J_n(x,r)$. 

Let $\alpha$ be the constant so that $(G,\rho_n)$ is $\alpha$-roughly geodesic. We show that each point $z \in J_n(x,r)$ is within $\alpha+2$ of a point in $B_n(x,r)$. Pick $z \in J_n(x,r)$, we will first show that $\rho_n(x,z) \leq r+2.$ If $\rho_n(x,z) > r+2$, there exists a chain $c = \{ h_1, \dots, h_{r+2} \}$ of length $r+2$ separating $x$ and $z$. By the definition of $\mathcal{H}^r(x)$, we have $h_{r+2} \in \mathcal{H}^r(x)$. Since $z \in h^+$ for each $h \in \mathcal{H}^r(x)$ with $x \in h^+$, it follows that $z \in h_{r+2}^+$, which contradicts the fact that $h_{r+2}$ separates $x$ and $z$. This shows that $\rho_n(x,z) \leq r+2$. Since $\sigma_n$ is $\alpha$-roughly geodesic, we have $\rho_n(z,w) \leq 2+\alpha$ for some $w \in B_n(x,r).$
\end{proof}

The above lemma shows that the balls $B_n(x,r)$ are $A$-convex for a uniform $A$ since they are at a uniform Hausdorff distance from the $E$-convex set $J_n(x,r)$. Now, Theorem J of \cite{haettelhodapetyt:coarse} implies that any collection of pairwise-intersecting balls in this metric must admit a point that is close to all of them. Since $(G, \rho_n)$ is roughly geodesic, this completes the proof of coarse injectivity for the metrics $(G, \rho_n)$. 

Passing to the injective hull of $(G, \rho_n)$ yields a geometric action of $G$ on an injective metric space for each $n$. It is worth noting that the points $(G, \rho_n)$ and the injective hull $E(G, \rho_n)$ represent the same point in $\DDGS$, since $(G, \rho_n)$ is coarsely dense in $E(G,\rho_n)$.

\subsection{Going to infinity in $\DDGS$ by strangling our favorite subsurface}

Having constructed infinitely many injective metric spaces—one for each $n$—we now aim to show that they yield an infinite diameter subset of $\DDGS$. Recall that when we defined the injective metric $(G, \rho_n)$, we chose intervals $I$ of length $10E n$ along geodesics $b \subset \Sigma$, and of length $10E$ along geodesics $b \subset Y$, where $Y$ is a proper subsurface. As $n \to \infty$, the thickness of the $\Sigma(n)$-curtains in $(G, \rho_1)$ grows unbounded; however, each of these $\Sigma(n)$-curtains contributes only 1 unit of distance in the metric $\rho_n$. 

Let $g$ be the pseudo-Anosov we fixed in Section \ref{sec:fixed_p.A}. Now, pick a point $p \in \mathcal{C}(\Sigma)$ so that $\beta:=\{g^ip\}$ is an $E$-rough-geodesic axis (for any $E$-geodesic hyperbolic space, and a loxodromic isometry $g$ on it, there exists a point $p$ such that $\{g^ip\}$ is an $E'$-rough geodesic where $E'$ depends only on $E$). The distance in $\mathcal{C}(\Sigma)$ between $p$ and $g^k.p$ is roughly $k$, and hence it is approximately $k$ in $(G, \rho_1)$ (by approximately, we allow errors depending on $E$ \emph{and} on $g$). However, when $G$ is equipped with the metric $\rho_n$, the distance $\rho_n(e, g^k)$ is approximately $\frac{k}{n}$ (indeed, by Definition \ref{def:weights}, for any $n,k \in \mathbb{N}$ if a chain $c$ realizes the distance $\rho_n(e,g^k)$, then every member of $c$ is labeled by $\mathcal{C}(\Sigma).$) In particular, this implies that 
\[
\frac{\ell_{\rho_1}(g)}{\ell_{\rho_n}(g)} \to \infty \quad \text{as} \quad n \to \infty.
\]

Now, recall that the distance on the space $\DDGS$ is given by
\[
D(\rho_1, \rho_n) := \log \sup_{g \in G_0} \frac{\ell_{\rho_1}(g)}{\ell_{\rho_n}(g)} + \log \sup_{g \in G_0} \frac{\ell_{\rho_n}(g)}{\ell_{\rho_1}(g)}.
\]
The above argument shows that the first term in the sum diverges to infinity. However, it is a priori possible that the second term approaches negative infinity. We now argue that this does not happen.

The key observation is that for an element $h$ that is active on a proper subsurface $Y$, we have $\ell_{\rho_i}(h) = \ell_{\rho_j}(h)$ for all $i \neq j$, since the width of curtains originating from a proper subsurface $Y$ remains unchanged when modifying $\rho_i$. Consequently, we obtain 
\[
\sup_{g \in G_0} \frac{\ell_{\rho_n}(g)}{\ell_{\rho_1}(g)} \geq 1.
\]
This proves Theorem~\ref{thm:infinitely_many_injective}.

    \bibliography{main}

\newcommand{\noop}[1]{}
\begin{thebibliography}{LBMB22}

\bibitem[ACT15]{Arzhantseva2015}
Goulnara Arzhantseva, Christopher Cashen, and Jing Tao.
\newblock Growth tight actions.
\newblock {\em Pacific Journal of Mathematics}, 278(1):1--49, September 2015.

\bibitem[AL02]{AL02:GrowthTight}
G.N Arzhantseva and I.G. Lysenok.
\newblock Growth tightness for word hyperbolic groups.
\newblock {\em Mathematische Zeitschrift}, 241(3):597--611, 2002.

\bibitem[Bow18]{bowditch:convex}
Brian~H. Bowditch.
\newblock Convex hulls in coarse median spaces.
\newblock {\em Preprint available at \mbox{homepages.warwick.ac.uk/\texttildelow masgak/papers/hulls-cms.pdf}}, 2018.

\bibitem[CGYZ24]{CGZ24:Confined}
Inhyeok Choi, Ilya Gekhtman, Wenyuan Yang, and Tianyi Zheng.
\newblock Confined subgroups in groups with contracting elements, 2024.

\bibitem[Cho21]{Choi2021PseudoAnosovsAE}
Inhyeok Choi.
\newblock Pseudo-anosovs are exponentially generic in mapping class groups.
\newblock {\em arXiv preprint arXiv:2110.06678}, 2021.

\bibitem[Cho24]{Full}
Inhyeok Choi.
\newblock Counting pseudo-anosovs as weakly contracting isometries, 2024.

\bibitem[CMGR]{CMGR:GreenDensity}
Stephen Cantrell, Didac Mart\'inez-Granado, and Eduardo Reyes.
\newblock Density of green metrics for hyperbolic groups.

\bibitem[COR22]{COR22:Manhattan}
Stephen Cantrell and Eduardo Oregón-Reyes.
\newblock Manhattan geodesics and the boundary of the space of metric structures on hyperbolic groups, 2022.

\bibitem[Cou24]{COULON2024}
R\'emi Coulon.
\newblock Patterson–sullivan theory for groups with a strongly contracting element.
\newblock {\em Ergodic Theory and Dynamical Systems}, page 1–56, March 2024.

\bibitem[CS15]{charneysultan:contracting}
Ruth Charney and Harold Sultan.
\newblock Contracting boundaries of {$\rm CAT(0)$} spaces.
\newblock {\em J. Topol.}, 8(1):93--117, 2015.

\bibitem[CT16]{CT16:GrowthTight_Word}
Christopher~H. Cashen and Jing Tao.
\newblock Growth tight actions of product groups.
\newblock {\em Groups, Geometry and Dynamics}, 10(2):753–770, 2016.

\bibitem[CW21]{calvezwiest:morse}
Matthieu Calvez and Bert Wiest.
\newblock Morse elements in {G}arside groups are strongly contracting.
\newblock {\em arXiv:2106.14826}, 2021.

\bibitem[DOP00]{DOP00:GrowthTight_Counter}
F.~Dal'bo, JP. Otal, and M.~Peigné.
\newblock Séries de poincaré des groupes géométriquement finis.
\newblock {\em Israel Journal of Mathematics}, 118(4):109--124, 2000.

\bibitem[DPPS11]{DPPS11:GrowthTight}
Françoise Dal’bo, Marc Peign\'e, Jean-Claude Picaud, and Andrea Sambusetti.
\newblock On the growth of quotients of kleinian groups.
\newblock {\em Ergodic Theory and Dynamical Systems}, 31(3):835–851, 2011.

\bibitem[DY24]{DY24:GrowthTight}
Lihuang Ding and Wenyuan Yang.
\newblock Growth tightness of quotients by confined subgroups, 2024.

\bibitem[Far06]{FarbProblems}
Benson Farb, editor.
\newblock {\em Problems on mapping class groups and related topics}, volume~74 of {\em Proceedings of Symposia in Pure Mathematics}.
\newblock American Mathematical Society, Providence, RI, 2006.

\bibitem[FG23]{FG23:ConfinedInfiniteVolume}
Mikolaj Fraczyk and Tsachik Gelander.
\newblock {Infinite volume and infinite injectivity radius}.
\newblock {\em Annals of Mathematics}, 197(1):389 -- 421, 2023.

\bibitem[Fur02]{Fur02:Coarse}
Alex Furman.
\newblock {\em Coarse-Geometric Perspective on Negatively Curved Manifolds and Groups}, pages 149--166.
\newblock Springer Berlin Heidelberg, Berlin, Heidelberg, 2002.

\bibitem[GdlH97]{growth_tight}
R.~Grigorchuk and P.~de~la Harpe.
\newblock On problems related to growth, entropy, and spectrum in group theory.
\newblock {\em J. Dynam. Control Systems}, 3(1):51--89, 1997.

\bibitem[GL24]{GL24:ConfinedStationary}
Ilya Gekhtman and Arie Levit.
\newblock Stationary random subgroups in negative curvature, 2024.

\bibitem[Gri77]{Gri77:Symmetric}
R.~Grigorchuk.
\newblock Symmetric random walks on discrete groups.
\newblock {\em Uspekhi Mat. Nauk}, pages 217--218, 1977.

\bibitem[GY22]{Gekhtman2022}
Ilya Gekhtman and Wen‐yuan Yang.
\newblock Counting conjugacy classes in groups with contracting elements.
\newblock {\em Journal of Topology}, 15(2):620–665, May 2022.

\bibitem[HHP23]{haettelhodapetyt:coarse}
Thomas Haettel, Nima Hoda, and Harry Petyt.
\newblock Coarse injectivity, hierarchical hyperbolicity and semihyperbolicity.
\newblock {\em Geom. Topol.}, 27(4):1587--1633, 2023.

\bibitem[HZ97]{HZ97:Confined}
B.~Hartley and A.~E. Zalesskiĭ.
\newblock Confined subgroups of simple locally finite groups and ideals of their group rings.
\newblock {\em Journal of the London Mathematical Society}, 55(2):210--230, 04 1997.

\bibitem[JJ20]{JM20:FreeGap}
Matsuzaki~Katsuhiko Jaerisch~Johannes.
\newblock Growth and cogrowth of normal subgroups of a free group.
\newblock {\em Groups Geometry and Dynamics}, 14(2):349--368, 2020.

\bibitem[Lan13a]{Lang13:Injective}
Urs Lang.
\newblock Injective hulls of certain discrete metric spaces and groups.
\newblock {\em J. Topol. Anal.}, 5(3):297--331, 2013.

\bibitem[Lan13b]{L13:Injective}
Urs Lang.
\newblock Injective hulls of certain discrete metric spaces and groups.
\newblock {\em J. Topol. Anal.}, 5(3):297--331, 2013.

\bibitem[LBMB18]{LM18:Confined2}
Adrien Le~Boudec and Nicolás Matte~Bon.
\newblock Subgroup dynamics and $c^\ast $-simplicity of groups of homeomorphisms.
\newblock {\em Annales scientifiques de l'École Normale Supérieure}, 51(3):557--602, 2018.

\bibitem[LBMB22]{LM22:Confined3}
Adrien Le~Boudec and Nicolás Matte~Bon.
\newblock Confined subgroups and high transitivity.
\newblock {\em Annales Henri Lebesgue}, 51(5):491--522, 2022.

\bibitem[LRT11]{Lenzhen2011}
Anna Lenzhen, Kasra Rafi, and Jing Tao.
\newblock Bounded combinatorics and the lipschitz metric on teichm\"{u}ller space.
\newblock {\em Geometriae Dedicata}, 159(1):353–371, October 2011.

\bibitem[Mah11]{maher:random}
Joseph Maher.
\newblock Random walks on the mapping class group.
\newblock {\em Duke Math. J.}, 156(3):429--468, 2011.

\bibitem[Min96]{Minsky1996}
Yair~N. Minsky.
\newblock Extremal length estimates and product regions in {T}eichmüller space.
\newblock {\em Duke Math. J.}, 83(2):249--286, 05 1996.

\bibitem[MM99]{masurminsky:geometry:1}
Howard~A. Masur and Yair~N. Minsky.
\newblock Geometry of the complex of curves.~{I}. {H}yperbolicity.
\newblock {\em Invent. Math.}, 138(1):103--149, 1999.

\bibitem[MM00]{masurminsky:geometry:2}
Howard~A. Masur and Yair~N. Minsky.
\newblock Geometry of the complex of curves.~{II}. {H}ierarchical structure.
\newblock {\em Geom. Funct. Anal.}, 10(4):902--974, 2000.

\bibitem[MYJ19]{JMY19:Normalizer}
Katsuhiko Matsuzaki, Yasuhiro Yabuki, and Johannes Jaerisch.
\newblock Normalizer, divergence type and patterson measure for discrete groups of the gromov hyperbolic space, 2019.

\bibitem[OR22]{OR22:SpaceMetric}
Eduardo Oregón-Reyes.
\newblock The space of metric structures on hyperbolic groups, 2022.

\bibitem[PSZ22]{PSZCAT}
Harry Petyt, Davide Spriano, and Abdul Zalloum.
\newblock Hyperbolic models for cat (0) spaces.
\newblock {\em arXiv preprint arXiv:2207.14127}, 2022.

\bibitem[PZ24]{PetytZalloum24}
Harry Petyt and Abdul Zalloum.
\newblock Constructing metric spaces from systems of walls.
\newblock {\em arXiv:2404.12057}, 2024.
\newblock With an appendix together with Davide Spriano.

\bibitem[Rey23]{Reyes23:Thesis}
Eduardo Reyes.
\newblock {\em Some aspects of geometric actions of hyperbolic and relatively hyperbolic groups}.
\newblock PhD thesis, UC Berkeley, 2023.

\bibitem[RST22]{Russell2018ConvexityIH}
Jacob Russell, Davide Spriano, and Hung~Cong Tran.
\newblock Convexity in hierarchically hyperbolic spaces.
\newblock {\em To appear in Algebr. Geom. Topol.}, 2022.

\bibitem[RV21]{Rafi2021}
Kasra Rafi and Yvon Verberne.
\newblock Geodesics in the mapping class group.
\newblock {\em Algebr. Geom. Topol.}, 21(6):2995--3017, 2021.

\bibitem[Sam02]{Sam02:GrowthTight_Surface}
Andrea Sambusetti.
\newblock Growth tightness of surface groups.
\newblock {\em Expositiones Mathematicae}, 20(4):345--363, 2002.

\bibitem[Sam08]{Sam08:GrowthTight_Cocompact}
Andrea Sambusetti.
\newblock {Asymptotic properties of coverings in negative curvature}.
\newblock {\em Geometry \& Topology}, 12(1):617 -- 637, 2008.

\bibitem[SZ]{soergelzalloum:morse}
Mireille Soergel and Abdul Zalloum.
\newblock Morse geodesics in weakly modular graphs.

\bibitem[SZ24]{SistoZalloum2024}
Alessandro Sisto and Abdul Zalloum.
\newblock Morse subsets of injective spaces are strongly contracting.
\newblock {\em Groups, Geometry, and Dynamics}, August 2024.

\bibitem[WXY24]{W24:Marked}
Renxing Wan, Xiaoyu Xu, and Wenyuan Yang.
\newblock Marked length spectrum rigidity in groups with contracting elements, 2024.

\bibitem[Yan14]{Yan14:GrowthTight}
Wen-Yuan Yang.
\newblock Growth tightness for groups with contracting elements.
\newblock {\em Mathematical Proceedings of the Cambridge Philosophical Society}, 157(2):297–319, 2014.

\bibitem[Yan18]{Yang17:ConvexCocompact}
Wen-yuan Yang.
\newblock Statistically convex-cocompact actions of groups with contracting elements, 02 2018.

\bibitem[Yan20]{Yang2018}
Wen-yuan Yang.
\newblock Genericity of contracting elements in groups.
\newblock {\em Math. Ann.}, 376(3-4):823--861, 2020.

\bibitem[Yan22]{Yang:PS&growth}
Wen-Yuan Yang.
\newblock Patterson–sullivan measures and growth of relatively hyperbolic groups.
\newblock {\em Peking Mathematical Journal}, 5:153--212, 2022.

\bibitem[Zal23]{zalloum:injectivity}
Abdul Zalloum.
\newblock Injectivity, cubical approximations and equivariant wall structures beyond {$\rm CAT(0)$} cube complexes.
\newblock {\em arXiv:2305.02951}, 2023.

\bibitem[Zbi23]{zbinden:small}
Stefanie Zbinden.
\newblock Small cancellation groups with and without sigma-compact {M}orse boundary.
\newblock {\em arXiv preprint arXiv:2307.13325}, 2023.

\end{thebibliography}
	\bibliographystyle{alpha}
\end{document}